\documentclass[12pt,a4paper,reqno]{amsart}

\usepackage[margin=1in,footskip=0.25in]{geometry}
\usepackage{amsmath, amsthm, amssymb}
\usepackage{graphicx}
\usepackage{mathrsfs}

\usepackage{float}

\usepackage{booktabs}

\usepackage{hyperref}
\usepackage{xcolor}

\newtheorem{theorem}{Theorem}

\date{}

\allowdisplaybreaks
\begin{document}

\title[Discontinuous piecewise differential systems with nilpotent saddles]
{Crossing limit cycles of discontinuous piecewise differential systems with nilpotent saddles separated by a nonregular line}

\author{Sonia Isabel Renteria Alva$^1$}
\address{$^1$ Instituto de Matemática Pura e Aplicada, Estrada Dona Castorina 110, Jardim Botânico, Rio de Janeiro, 22460-320, Brazil}
\email{sonia.alva@impa.br}

\author{Pedro Iván Suárez Navarro$^2$}
\address{$^2$ Instituto de Matemática Pura e Aplicada,  Estrada Dona Castorina 110, Jardim Botânico, Rio de Janeiro, 22460-320, Brazil}
\email{ivan.suarez@impa.br}



\begin{abstract}
In this paper, we study crossing limit cycles of planar discontinuous piecewise differential systems separated by a nonregular switching line, where one subsystem is a linear differential center and the other belongs to one of six families of Hamiltonian nilpotent saddle systems. This work extends previous results for discontinuous piecewise differential systems separated by a straight line to the more general setting of nonregular switching boundaries. The problem is closely related to the second part of Hilbert’s 16th problem for discontinuous piecewise differential systems. We prove that the systems under consideration admit at most seven crossing limit cycles. Moreover, for each family of Hamiltonian nilpotent saddle systems, we provide explicit examples exhibiting four crossing limit cycles.
\end{abstract}

\maketitle

\section{Introduction and statement of the main results}

Planar discontinuous piecewise differential systems have become an active and
rapidly developing area of research in recent years. This growing interest is
largely motivated by their close relationship with extensions of the second part
of Hilbert’s 16th problem, as well as by their ability to describe a wide
range of phenomena arising in real-world applications. Comprehensive treatments
of the theory and applications of such systems can be found, for instance, in
the monographs~\cite{bernardo2008piecewise,simpson2010bifurcations}.

From a qualitative point of view, the analysis of piecewise systems presents substantial challenges, even when the discontinuity set is given by simple geometric objects such as straight lines, circles, or low-degree algebraic curves. The presence of different vector fields acting on either side of the
discontinuity may lead to complex dynamical behaviors that do not occur in smooth
systems. In particular, the interaction of the vector fields across the discontinuity manifold plays a crucial role in the appearance and bifurcation of
periodic solutions; see, for example,
\cite{anacleto2021limit,llibre2022crossing,alva2025crossing,baymout2024limit}.

In this work, we extend the results obtained in~\cite{benabdallah2022four} for discontinuous piecewise differential systems separated by a straight line to the case in which the discontinuity curve is a nonregular line. Systems separated by such curves have recently attracted attention, and related contributions under different classes of vector fields can be found in
\cite{cabrera2026limit,baymout2024limit,he2024limit,huan2019limit,berbache2025lower,he2025limit, llibre2025algebraic, chen2025limit,chen2025establishing}.

A key ingredient in our analysis is the normal form of linear differential systems possessing a center. As shown in~\cite{llibre2018piecewise}, any such system can be written in the form
\begin{equation}\label{sistemalinear}
(\mathtt{L}_c):\;
\begin{aligned}
\dot{x} &= -(A^2+\omega^2)y - A x + B,\dot{y} &= A y + x + C,\quad \text{ with } \omega>0 
\end{aligned}
\end{equation}
with the first integral
\begin{equation}\label{intprimeiralinear}
H_L(x,y) = (A y+x)^2+2(C x -B y)+\omega^2 y^2.
\end{equation}

Colak, Llibre, and Vals~\cite{colak2014hamiltonian} provided normal forms and global
phase portraits in the Poincar\'e disk for all Hamiltonian nilpotent centers of
linear plus cubic homogeneous planar polynomial vector fields. Among these
systems, six classes of Hamiltonian nilpotent saddle vector fields play a central
role in our setting. For completeness, we list below these six classes together
with their corresponding first integrals.

\begin{description}
\item[$(\mathtt{N}_1)$]
$\dot{x}=ax+by,\ \dot{y}=-\frac{a^{2}}{b}x-ay+x^{3},\quad b<0$,
with the first integral
$$\widetilde{H}_1(x,y)=-\frac{x^{4}}{4}
+\frac{a}{2b}x^{2}
+\frac{b}{2}y^{2}
+a x y.$$
\item[$(\mathtt{N}_2)$]
$\dot{x}=ax+by-x^{3},\ \dot{y}=-\frac{a^{2}}{b}x-ay+3x^{2}y,\quad a>0$,
with the first integral
$$\widetilde{H}_2(x,y)=- x^{3}y
+\frac{a^{2}}{2b}x^{2}
+\frac{b}{2}y^{2}
+a x y.$$
\item[$(\mathtt{N}_3)$]
$\dot{x}=ax+by-3x^{2}y+y^{3},\ 
 \dot{y}=\left(c-\frac{a^{2}}{b+c}\right)x-ay+3xy^{2},$
either $a=b=0,\ c<0$, with the first integral
$$\widetilde{H}_3^{1}(x,y)=\frac{y^{4}}{4}
-\frac{3}{2}x^{2}y^{2}
-\frac{c}{2}x^{2},$$
or $c=0,\ ab\neq0,\ \frac{a^{2}}{b}-6b>0$,  with the first integral
$$\widetilde{H}_3^{2}(x,y)=\frac{y^{4}}{4}
-\frac{3}{2}x^{2}y^{2}
+\frac{a^{2}}{2b}x^{2}
+\frac{b}{2}y^{2}
+a x y.$$
\item[$(\mathtt{N}_4)$]
$\dot{x}=ax+by-3x^{2}y-y^{3},\ 
 \dot{y}=\left(c-\frac{a^{2}}{b+c}\right)x-ay+3xy^{2},$
either $a=b=0,\ c>0$, with the first integral 
$$ \widetilde{H}_4^{1}(x,y)=-\frac{y^{4}}{4}
-\frac{3}{2}x^{2}y^{2}
-\frac{c}{2}x^{2},$$ or $c=0,\ a\neq0,\ b<0$, with the first integral
$$ \widetilde{H}_4^{2}(x,y)=-\frac{y^{4}}{4}
-\frac{3}{2}x^{2}y^{2}
+\frac{a^{2}}{2b}x^{2}
+\frac{b}{2}y^{2}
+a x y.$$
\item[$(\mathtt{N}_5)$]
$\dot{x}=ax+by-3\mu x^{2}y+y^{3},\ 
 \dot{y}=\left(c-\frac{a^{2}}{b+c}\right)x-ay+x^{3}+3\mu xy^{2},$
either $a=b=0,\ c<0$, with the first integral 
$$ \widetilde{H}_5^{1}(x,y)=\frac{y^{4}-x^{4}}{4}
-\frac{3}{2}\mu x^{2}y^{2}
-\frac{c}{2}x^{2},$$ 
or $c=0,\ b\neq0,\
\dfrac{a^{4}-b^{4}-6a^{2}b^{2}\mu}{b}>0$, with the first integral
$$\widetilde{H}_5^{2}(x,y)=\frac{y^{4}-x^{4}}{4}
-\frac{3}{2}\mu x^{2}y^{2}
+\frac{a^{2}}{2b}x^{2}
+\frac{b}{2}y^{2}
+a x y.$$

\item[$(\mathtt{N}_6)$]
$\dot{x}=ax+by-3\mu x^{2}y-y^{3},\ 
 \dot{y}=\left(c-\frac{a^{2}}{b+c}\right)x-ay+x^{3}+3\mu xy^{2},$
either $a=b=0,\ c>0$, with the first integral 
$$\widetilde{H}_6^{1}(x,y)=-\frac{x^{4}+y^{4}}{4}
-\frac{3}{2}\mu x^{2}y^{2}
-\frac{c}{2}x^{2},$$ or $c=0,\ b\neq0,\
\dfrac{a^{4}+b^{4}+6a^{2}b^{2}\mu}{b}<0$, with the first integral
$$\widetilde{H}_6^{2}(x,y)=-\frac{x^{4}+y^{4}}{4}
-\frac{3}{2}\mu x^{2}y^{2}
+\frac{a^{2}}{2b}x^{2}
+\frac{b}{2}y^{2}
+a x y.$$
\end{description}

The main objective of this paper is to determine the maximum number of crossing
limit cycles that can occur in a class of discontinuous piecewise differential
systems separated by the nonregular line
$\Sigma=\Sigma^{+}\cap \Sigma^{-} $. 
In this class, the dynamics are governed by
a linear differential center $(\mathtt{L}_c)$ in the region
$\Sigma^-=\{(x,y)\in \mathbb{R}^2: x\leq 0\} \cup \{(x,y)\in \mathbb{R}^2: x\geq 0, y\leq 0\}$, while in the
region $\Sigma^+=\{(x,y) \in \mathbb{R}^2: x\geq 0, y\geq 0\}$, the system is given by
one of the six Hamiltonian nilpotent saddle systems
$(\mathtt{N}_1),\dots,(\mathtt{N}_6)$.

Our main result establishes a uniform upper bound for the number of crossing limit
cycles arising in such systems and shows that this bound is sharp for each of the
six classes considered.

\begin{theorem}\label{Teor-Principal}
Consider discontinuous piecewise differential systems separated by a nonregular line $\Sigma$,
formed by a linear differential center \eqref{sistemalinear}  and by one
of the six classes of Hamiltonian nilpotent saddles $(\mathtt{N}_k)$,
$k\in\{1,\dots,6\}$.
Then such systems admit at most seven crossing limit cycles.
Moreover, for each class there exist examples exhibiting  four crossing limit cycles
(see Figures~\ref{fig-Hc-H1}--\ref{fig-Hc-H62}).
\end{theorem}



Table~\ref{TableLCy} summarizes the results obtained in
Theorem~\ref{Teor-Principal}. For each of the six classes of discontinuous
piecewise differential systems considered in this work, the table reports the
corresponding upper bound for the number of crossing limit cycles. The numbers
shown in parentheses indicate the maximal number of limit cycles that can
actually be achieved by suitable choices of parameters. For $k=3,4,5,6$, the
labels $N_k^{1}$ and $N_k^{2}$ distinguish two subcases within the class $N_k$.
Explicit constructions realizing these maximal configurations are given in the
proof of Theorem~\ref{Teor-Principal}.

\begin{table}[H]

\centering
\begin{tabular}{l | c c c c c c c c c r}
  & $\mathtt{N}_1$  & $\mathtt{N}_2$   & $\mathtt{N}_3^{1}$  &  $\mathtt{N}_3^{2}$  &  $\mathtt{N}_4^{1}$  &  $\mathtt{N}_4^{2}$ &  $\mathtt{N}_5^{1}$ &  $\mathtt{N}_5^{2}$ &  $\mathtt{N}_6^{1}$ &  $\mathtt{N}_6^{2}$\\
  \hline
  $\mathtt{L}_c$  & 7(4) & 7(4)  & 7(4) & 7(4) & 7(4) & 7(4) & 7(4) & 7(4) & 7(4) & 7(4)\\
\end{tabular}
\caption{Upper bounds for the number of crossing limit cycles obtained in
Theorem~\ref{Teor-Principal}.}
\label{TableLCy}
\end{table}


In Section \ref{sect:02} we present the Hamiltonian nilpotent saddles after an affine
change of variables. Theorem \ref{Teor-Principal} is proved in Section \ref{sect:03}.

\section{The Hamiltonian nilpotent saddles ($\mathtt{N}_1$), ($\mathtt{N}_2$), ($\mathtt{N}_3$),($\mathtt{N}_4$),($\mathtt{N}_5$) and ($\mathtt{N}_6$) after an affine change of variables}
\label{sect:02}
In this section, we give the expression of the Hamiltonian nilpotent saddles ($\mathtt{N}_1$), ($\mathtt{N}_2$), ($\mathtt{N}_3$),($\mathtt{N}_4$),($\mathtt{N}_5$) and ($\mathtt{N}_6$), as well as to their first integrals after the general affine change of
variables $(x, y) \rightarrow (a_1 x +b_1 y + c_1, \alpha_1 x + \beta_1 y + \gamma_1)$,
with $b_1 \alpha_1 - a_1 \beta_1\neq 0$.
Thus, after this affine change of variables the differential system ($\mathtt{N}_1$) becomes 
\begin{equation}\label{eq:system}
\begin{aligned}
\dot{x} &=
\frac{1}{b\left(b_{1}\alpha_{1}-a_{1}\beta_{1}\right)}
\big(
-a^{2} b_{1} c_{1} + b b_{1} c_{1}^{3}- a b b_{1}\gamma_{1}
- b^{2}\beta_{1}\gamma_{1} - a b c_{1}\beta_{1}  \\
&\quad
+ 3 a_{1}^{2} b b_{1} c_{1} x^{2}
+ a_{1}^{3} b b_{1} x^{3}
+ 6 a_{1} b b_{1}^{2} c_{1} x y
+ 3 a_{1}^{2} b b_{1}^{2} x^{2} y  \\
&\quad
+ 3 b b_{1}^{3} c_{1} y^{2}
+ 3 a_{1} b b_{1}^{3} x y^{2}
+ b b_{1}^{4} y^{3}  \\
&\quad
+ \left(
-a^{2} a_{1} b_{1}
+ 3 a_{1} b b_{1} c_{1}^{2}
- a b b_{1}\alpha_{1}
- a a_{1} b\beta_{1}
- b^{2}\alpha_{1}\beta_{1}
\right)x  \\
&\quad
+ \left(
-a^{2} b_{1}^{2}
+ 3 b b_{1}^{2} c_{1}^{2}
- 2 a b b_{1}\beta_{1}
- b^{2}\beta_{1}^{2}
\right)y
\big), \\[1ex]
\dot{y} &=
\frac{1}{b\left(b_{1}\alpha_{1}-a_{1}\beta_{1}\right)}
\big(
a^{2} a_{1} c_{1} - a_{1} b c_{1}^{3} + a a_{1} b\gamma_{1}
+ b^{2}\alpha_{1}\gamma_{1} + a b c_{1}\alpha_{1}  \\
&\quad
- 3 a_{1}^{3} b c_{1} x^{2}
- a_{1}^{4} b x^{3}
- 6 a_{1}^{2} b b_{1} c_{1} x y
- 3 a_{1}^{3} b b_{1} x^{2} y  \\
&\quad
- 3 a_{1} b b_{1}^{2} c_{1} y^{2}
- 3 a_{1}^{2} b b_{1}^{2} x y^{2}
- a_{1} b b_{1}^{3} y^{3}  \\
&\quad
+ \left(
a^{2} a_{1}^{2}
- 3 a_{1}^{2} b c_{1}^{2}
+ 2 a a_{1} b\alpha_{1}
+ b^{2}\alpha_{1}^{2}
\right)x  \\
&\quad
+ \left(
a^{2} a_{1} b_{1}
- 3 a_{1} b b_{1} c_{1}^{2}
+ a b b_{1}\alpha_{1}
+ a a_{1} b\beta_{1}
+ b^{2}\alpha_{1}\beta_{1}
\right)y
\big),
\end{aligned}
\end{equation}


with the first integral
\begin{equation}\label{eq:Htilde1}
\begin{aligned}
H_1(x,y)
&=\frac{a}{2b}\left(c_{1}+a_{1}x+b_{1}y\right)^{2}
-\frac{1}{4}\left(c_{1}+a_{1}x+b_{1}y\right)^{4}  \\
&\quad
+a\left(c_{1}+a_{1}x+b_{1}y\right)
\left(x\alpha_{1}+y\beta_{1}+\gamma_{1}\right)
+\frac{b}{2}\left(x\alpha_{1}+y\beta_{1}+\gamma_{1}\right)^{2}.
\end{aligned}
\end{equation}

The differential system ($\mathtt{N}_2$) becomes
\begin{equation}\label{eq:systemC2}
\begin{aligned}
\dot{x} &=
\frac{1}{b\left(b_{1}\alpha_{1}-a_{1}\beta_{1}\right)}
\big(
-a^{2} b_{1} c_{1} - a b c_{1}\beta_{1} + b c_{1}^{3}\beta_{1}
- a b b_{1}\gamma_{1} + 3 b b_{1} c_{1}^{2}\gamma_{1}
- b^{2}\beta_{1}\gamma_{1}  \\
&\quad
+ 4 b b_{1}^{3}\beta_{1} y^{3}
+ a_{1}^{2} b\left(3 b_{1}\alpha_{1}+a_{1}\beta_{1}\right) x^{3}
+ 6\left(a_{1} b b_{1}^{2}\alpha_{1}
+ a_{1}^{2} b b_{1}\beta_{1}\right)x^{2}y  \\
&\quad
+ 3\left(b b_{1}^{3}\alpha_{1}
+ 3 a_{1} b b_{1}^{2}\beta_{1}\right)x y^{2}
+ 3 a_{1} b\left(
2 b_{1} c_{1}\alpha_{1}
+ a_{1} c_{1}\beta_{1}
+ a_{1} b_{1}\gamma_{1}
\right)x^{2}  \\
&\quad
+ 6 b b_{1}\left(
b_{1} c_{1}\alpha_{1}
+ 2 a_{1} c_{1}\beta_{1}
+ a_{1} b_{1}\gamma_{1}
\right)x y
+ 3\left(
3 b b_{1}^{2} c_{1}\beta_{1}
+ b b_{1}^{3}\gamma_{1}
\right)y^{2}  \\
&\quad
+ \big(
-a^{2} a_{1} b_{1}
- a b b_{1}\alpha_{1}
+ 3 b b_{1} c_{1}^{2}\alpha_{1}
- a a_{1} b\beta_{1}
+ 3 a_{1} b c_{1}^{2}\beta_{1}
- b^{2}\alpha_{1}\beta_{1} \\
&\quad
+ 6 a_{1} b b_{1} c_{1}\gamma_{1}
\big)x 
+ \big(
-a^{2} b_{1}^{2}
- 2 a b b_{1}\beta_{1}
+ 6 b b_{1} c_{1}^{2}\beta_{1}
- b^{2}\beta_{1}^{2}
+ 6 b b_{1}^{2} c_{1}\gamma_{1}
\big)y
\big), 
\end{aligned}
\end{equation}

\begin{equation}\label{eq:systemC2}
\begin{aligned}
\dot{y} &=
\frac{1}{b\left(b_{1}\alpha_{1}-a_{1}\beta_{1}\right)}
\big(
a^{2} a_{1} c_{1} + a b c_{1}\alpha_{1}- b c_{1}^{3}\alpha_{1}
+ a a_{1} b\gamma_{1}- 3 a_{1} b c_{1}^{2}\gamma_{1}
+ b^{2}\alpha_{1}\gamma_{1}  \\
&\quad
- 4 a_{1}^{3} b\alpha_{1} x^{3}
- 3\left(3 a_{1}^{2} b b_{1}\alpha_{1}
+ a_{1}^{3} b\beta_{1}\right)x^{2}y
- 6 a_{1} b b_{1}\left(b_{1}\alpha_{1}
+ a_{1}\beta_{1}\right)x y^{2}  \\
&\quad
- \left(b b_{1}^{3}\alpha_{1}
+ 3 a_{1} b b_{1}^{2}\beta_{1}\right)y^{3}
- 3\left(3 a_{1}^{2} b c_{1}\alpha_{1}
+ a_{1}^{3} b\gamma_{1}\right)x^{2}  \\
&\quad
- 6 a_{1} b\left(
2 b_{1} c_{1}\alpha_{1}
+ a_{1} c_{1}\beta_{1}
+ a_{1} b_{1}\gamma_{1}
\right)x y
- 3 b b_{1}(
b_{1} c_{1}\alpha_{1}
+ 2 a_{1} c_{1}\beta_{1} \\
&\quad
+ a_{1} b_{1}\gamma_{1}
)y^{2} 
+ \big(
a^{2} a_{1}^{2}
+ 2 a a_{1} b\alpha_{1}
- 6 a_{1} b c_{1}^{2}\alpha_{1}
+ b^{2}\alpha_{1}^{2}
- 6 a_{1}^{2} b c_{1}\gamma_{1}
\big)x  \\
&\quad
+ \big(
a^{2} a_{1} b_{1}
+ a b b_{1}\alpha_{1}
- 3 b b_{1} c_{1}^{2}\alpha_{1}
+ a a_{1} b\beta_{1}
- 3 a_{1} b c_{1}^{2}\beta_{1}
+ b^{2}\alpha_{1}\beta_{1}  \\
&\quad
- 6 a_{1} b b_{1} c_{1}\gamma_{1}
\big)y
\big),
\end{aligned}
\end{equation}
with the first integral
\begin{align}
 H_2(x,y)= &\frac{a^{2}}{2b}\left(c_{1}+a_{1}x+b_{1}y\right)^{2}
+a\left(c_{1}+a_{1}x+b_{1}y\right)
\left(x\alpha_{1}+y\beta_{1}+\gamma_{1}\right) \nonumber\\
&
-\left(c_{1}+a_{1}x+b_{1}y\right)^{3}
\left(x\alpha_{1}+y\beta_{1}+\gamma_{1}\right)
+\frac{b}{2}\left(x\alpha_{1}+y\beta_{1}+\gamma_{1}\right)^{2}.
\end{align}

The differential system ($\mathtt{N}_3^{1}$) becomes
\begin{equation}\label{eq:systemC3}
\begin{aligned}
\dot{x} &=
\frac{1}{-b_{1}\alpha_{1}+a_{1}\beta_{1}}
\big(
- b_{1} c c_{1}
- 3 c_{1}^{2}\beta_{1}\gamma_{1}
- 3 b_{1} c_{1}\gamma_{1}^{2}
+ \beta_{1}\gamma_{1}^{3}  
+ (
-3 a_{1} b_{1}\alpha_{1}^{2}\\
&\quad
- 3 a_{1}^{2}\alpha_{1}\beta_{1}
+ \alpha_{1}^{3}\beta_{1}
)x^{3}
- 3(
b_{1}^{2}\alpha_{1}^{2}
+ 4 a_{1} b_{1}\alpha_{1}\beta_{1}
+ a_{1}^{2}\beta_{1}^{2}
- \alpha_{1}^{2}\beta_{1}^{2}
)x^{2}y  \\
&\quad
- 3\left(
3 b_{1}^{2}\alpha_{1}\beta_{1}
+ 3 a_{1} b_{1}\beta_{1}^{2}
- \alpha_{1}\beta_{1}^{3}
\right)x y^{2}
+ \left(
-6 b_{1}^{2}\beta_{1}^{2}
+ \beta_{1}^{4}
\right)y^{3}  \\
&\quad
- 3\left(
b_{1} c_{1}\alpha_{1}^{2}
+ 2 a_{1} c_{1}\alpha_{1}\beta_{1}
+ 2 a_{1} b_{1}\alpha_{1}\gamma_{1}
+ a_{1}^{2}\beta_{1}\gamma_{1}
- \alpha_{1}^{2}\beta_{1}\gamma_{1}
\right)x^{2}  \\
&\quad
- 6\left(
2 b_{1} c_{1}\alpha_{1}\beta_{1}
+ a_{1} c_{1}\beta_{1}^{2}
+ b_{1}^{2}\alpha_{1}\gamma_{1}
+ 2 a_{1} b_{1}\beta_{1}\gamma_{1}
- \alpha_{1}\beta_{1}^{2}\gamma_{1}
\right)x y  \\
&\quad
- 3\left(
3 b_{1} c_{1}\beta_{1}^{2}
+ 3 b_{1}^{2}\beta_{1}\gamma_{1}
- \beta_{1}^{3}\gamma_{1}
\right)y^{2}  
+ \big(
- a_{1} b_{1} c
- 3 c_{1}^{2}\alpha_{1}\beta_{1}\\
&\quad
- 6 b_{1} c_{1}\alpha_{1}\gamma_{1}
- 6 a_{1} c_{1}\beta_{1}\gamma_{1}
- 3 a_{1} b_{1}\gamma_{1}^{2}
+ 3 \alpha_{1}\beta_{1}\gamma_{1}^{2}
\big)x  \\
&\quad
+ \big(
- b_{1}^{2} c
- 3 c_{1}^{2}\beta_{1}^{2}
- 12 b_{1} c_{1}\beta_{1}\gamma_{1}
- 3 b_{1}^{2}\gamma_{1}^{2}
+ 3 \beta_{1}^{2}\gamma_{1}^{2}
\big)y
\big),
\end{aligned}
\end{equation}


\begin{equation}\label{eq:ydotC3}
\begin{aligned}
\dot{y}=&
\frac{1}{-b_{1}\alpha_{1}+a_{1}\beta_{1}}
\big(a_{1} c c_{1}- x^{3}\alpha_{1}^{2}\big(-6 a_{1}^{2}+\alpha_{1}^{2}\big) + 3 x^{2} y\big(
3 a_{1} b_{1}\alpha_{1}^{2}
+ 3 a_{1}^{2}\alpha_{1}\beta_{1} \\
&\quad
- \alpha_{1}^{3}\beta_{1}
\big)
- y^{3}\beta_{1}\big(
-3 b_{1}^{2}\alpha_{1}
- 3 a_{1} b_{1}\beta_{1}
+ \alpha_{1}\beta_{1}^{2}
\big) 
+ 3 x y^{2}\big(
b_{1}^{2}\alpha_{1}^{2}  \\
&\quad
+ 4 a_{1} b_{1}\alpha_{1}\beta_{1}
+ a_{1}^{2}\beta_{1}^{2}
- \alpha_{1}^{2}\beta_{1}^{2}
\big)
+ 3 c_{1}^{2}\alpha_{1}\gamma_{1}
+ 3 a_{1} c_{1}\gamma_{1}^{2}
- \alpha_{1}\gamma_{1}^{3}  \\
&\quad
+ 3 x^{2}\big(
3 a_{1} c_{1}\alpha_{1}^{2}
+ 3 a_{1}^{2}\alpha_{1}\gamma_{1}
- \alpha_{1}^{3}\gamma_{1}
\big)
+ 6 x y\big(
b_{1} c_{1}\alpha_{1}^{2}
+ 2 a_{1} c_{1}\alpha_{1}\beta_{1} \\
&\quad
+ 2 a_{1} b_{1}\alpha_{1}\gamma_{1}
+ a_{1}^{2}\beta_{1}\gamma_{1}
- \alpha_{1}^{2}\beta_{1}\gamma_{1}
\big)  \\
&\quad
+ 3 y^{2}\big(
2 b_{1} c_{1}\alpha_{1}\beta_{1}
+ a_{1} c_{1}\beta_{1}^{2}
+ b_{1}^{2}\alpha_{1}\gamma_{1}
+ 2 a_{1} b_{1}\beta_{1}\gamma_{1}
- \alpha_{1}\beta_{1}^{2}\gamma_{1}
\big)  \\
&\quad
+ x\big(
a_{1}^{2} c
+ 3 c_{1}^{2}\alpha_{1}^{2}
+ 12 a_{1} c_{1}\alpha_{1}\gamma_{1}
+ 3 a_{1}^{2}\gamma_{1}^{2}
- 3 \alpha_{1}^{2}\gamma_{1}^{2}
\big)  \\
&\quad
+ y\big(
a_{1} b_{1} c
+ 3 c_{1}^{2}\alpha_{1}\beta_{1}
+ 6 b_{1} c_{1}\alpha_{1}\gamma_{1}
+ 6 a_{1} c_{1}\beta_{1}\gamma_{1}
+ 3 a_{1} b_{1}\gamma_{1}^{2}
- 3 \alpha_{1}\beta_{1}\gamma_{1}^{2}
\big)
\big),
\end{aligned}
\end{equation}

with the first integral
\begin{align}
H_3^{1}(x,y)&=-\frac{1}{2}c\left(c_{1}+a_{1}x+b_{1}y\right)^{2}
-\frac{3}{2}\left(c_{1}+a_{1}x+b_{1}y\right)^{2}
\left(x\alpha_{1}+y\beta_{1}+\gamma_{1}\right)^{2} \nonumber\\
&\quad
+\frac{1}{4}\left(x\alpha_{1}+y\beta_{1}+\gamma_{1}\right)^{4}.
\end{align}
The differential system ($\mathtt{N}_3^{2}$) becomes
\begin{equation}\label{eq:xdotC4}
\begin{aligned}
\dot{x} &=
\frac{1}{b\left(b_{1}\alpha_{1}-a_{1}\beta_{1}\right)}
\big(
- a^{2} b_{1} c_{1}
- a b c_{1}\beta_{1}
- a b b_{1}\gamma_{1}
- b^{2}\beta_{1}\gamma_{1} 
+ 3 b c_{1}^{2}\beta_{1}\gamma_{1}  \\
&\quad
+ 3 b b_{1} c_{1}\gamma_{1}^{2}
- b\beta_{1}\gamma_{1}^{3}
- b\alpha_{1}\big(
-3 a_{1} b_{1}\alpha_{1}
- 3 a_{1}^{2}\beta_{1}
+ \alpha_{1}^{2}\beta_{1}
\big)x^{3}  \\
&\quad
- b\beta_{1}^{2}\big(
-6 b_{1}^{2} + \beta_{1}^{2}
\big)y^{3}
+ 3\big(
b b_{1}^{2}\alpha_{1}^{2}
+ 4 a_{1} b b_{1}\alpha_{1}\beta_{1}
+ a_{1}^{2} b\beta_{1}^{2}
- b\alpha_{1}^{2}\beta_{1}^{2}
\big)x^{2}y  \\
&\quad
+ 3\big(
3 b b_{1}^{2}\alpha_{1}\beta_{1}
+ 3 a_{1} b b_{1}\beta_{1}^{2}
- b\alpha_{1}\beta_{1}^{3}
\big)x y^{2} 
+ 3\big(
b b_{1} c_{1}\alpha_{1}^{2}
+ 2 a_{1} b c_{1}\alpha_{1}\beta_{1} \\
&\quad
+ 2 a_{1} b b_{1}\alpha_{1}\gamma_{1}
+ a_{1}^{2} b\beta_{1}\gamma_{1}
- b\alpha_{1}^{2}\beta_{1}\gamma_{1}
\big)x^{2} 
+ 6\big(
2 b b_{1} c_{1}\alpha_{1}\beta_{1}
+ a_{1} b c_{1}\beta_{1}^{2}  \\
&\quad
+ b b_{1}^{2}\alpha_{1}\gamma_{1}
+ 2 a_{1} b b_{1}\beta_{1}\gamma_{1}
- b\alpha_{1}\beta_{1}^{2}\gamma_{1}
\big)x y 
+ 3\big(
3 b b_{1} c_{1}\beta_{1}^{2}
+ 3 b b_{1}^{2}\beta_{1}\gamma_{1}
- b\beta_{1}^{3}\gamma_{1}
\big)y^{2}  \\
&\quad
+ \big(
- a^{2} a_{1} b_{1}
- a b b_{1}\alpha_{1}
- a a_{1} b\beta_{1}
- b^{2}\alpha_{1}\beta_{1}
+ 3 b c_{1}^{2}\alpha_{1}\beta_{1}
+ 6 b b_{1} c_{1}\alpha_{1}\gamma_{1}  \\
&\quad
+ 6 a_{1} b c_{1}\beta_{1}\gamma_{1}
+ 3 a_{1} b b_{1}\gamma_{1}^{2}
- 3 b\alpha_{1}\beta_{1}\gamma_{1}^{2}
\big)x 
+ \big(
- a^{2} b_{1}^{2}
- 2 a b b_{1}\beta_{1}  
- b^{2}\beta_{1}^{2} \\
&\quad
+ 3 b c_{1}^{2}\beta_{1}^{2}
+ 12 b b_{1} c_{1}\beta_{1}\gamma_{1}
+ 3 b b_{1}^{2}\gamma_{1}^{2}
- 3 b\beta_{1}^{2}\gamma_{1}^{2}
\big)y
\big),
\end{aligned}
\end{equation}

\begin{equation}\label{eq:ydotC4}
\begin{aligned}
\dot{y} &=
\frac{1}{b\left(b_{1}\alpha_{1}-a_{1}\beta_{1}\right)}
\big(
a^{2} a_{1} c_{1}
+ a b c_{1}\alpha_{1}
+ a a_{1} b\gamma_{1}
+ b^{2}\alpha_{1}\gamma_{1}  \\
&\quad
- 3 b c_{1}^{2}\alpha_{1}\gamma_{1}
- 3 a_{1} b c_{1}\gamma_{1}^{2}
+ b\alpha_{1}\gamma_{1}^{3}
+ \big(
-6 a_{1}^{2} b\alpha_{1}^{2}
+ b\alpha_{1}^{4}
\big)x^{3}  \\
&\quad
- 3\big(
3 a_{1} b b_{1}\alpha_{1}^{2}
+ 3 a_{1}^{2} b\alpha_{1}\beta_{1}
- b\alpha_{1}^{3}\beta_{1}
\big)x^{2}y  \\
&\quad
- 3\big(
b b_{1}^{2}\alpha_{1}^{2}
+ 4 a_{1} b b_{1}\alpha_{1}\beta_{1}
+ a_{1}^{2} b\beta_{1}^{2}
- b\alpha_{1}^{2}\beta_{1}^{2}
\big)x y^{2}  \\
&\quad
+ \big(
-3 b b_{1}^{2}\alpha_{1}\beta_{1}
- 3 a_{1} b b_{1}\beta_{1}^{2}
+ b\alpha_{1}\beta_{1}^{3}
\big)y^{3}  \\
&\quad
- 3\big(
3 a_{1} b c_{1}\alpha_{1}^{2}
+ 3 a_{1}^{2} b\alpha_{1}\gamma_{1}
- b\alpha_{1}^{3}\gamma_{1}
\big)x^{2}  \\
&\quad
- 6\big(
b b_{1} c_{1}\alpha_{1}^{2}
+ 2 a_{1} b c_{1}\alpha_{1}\beta_{1}
+ 2 a_{1} b b_{1}\alpha_{1}\gamma_{1}
+ a_{1}^{2} b\beta_{1}\gamma_{1}
- b\alpha_{1}^{2}\beta_{1}\gamma_{1}
\big)x y  \\
&\quad
- 3\big(
2 b b_{1} c_{1}\alpha_{1}\beta_{1}
+ a_{1} b c_{1}\beta_{1}^{2}
+ b b_{1}^{2}\alpha_{1}\gamma_{1}
+ 2 a_{1} b b_{1}\beta_{1}\gamma_{1}
- b\alpha_{1}\beta_{1}^{2}\gamma_{1}
\big)y^{2}  \\
&\quad
+ \big(
a^{2} a_{1}^{2}
+ 2 a a_{1} b\alpha_{1}
+ b^{2}\alpha_{1}^{2}
- 3 b c_{1}^{2}\alpha_{1}^{2}
- 12 a_{1} b c_{1}\alpha_{1}\gamma_{1}
- 3 a_{1}^{2} b\gamma_{1}^{2}
+ 3 b\alpha_{1}^{2}\gamma_{1}^{2}
\big)x  \\
&\quad
+ \big(
a^{2} a_{1} b_{1}
+ a b b_{1}\alpha_{1}
+ a a_{1} b\beta_{1}
+ b^{2}\alpha_{1}\beta_{1}
- 3 b c_{1}^{2}\alpha_{1}\beta_{1}
- 6 b b_{1} c_{1}\alpha_{1}\gamma_{1}  \\
&\quad\quad
- 6 a_{1} b c_{1}\beta_{1}\gamma_{1}
- 3 a_{1} b b_{1}\gamma_{1}^{2}
+ 3 b\alpha_{1}\beta_{1}\gamma_{1}^{2}
\big)y
\big),
\end{aligned}
\end{equation}

with the first integral
\begin{align}
H_3^{2}(x,y)
&= \frac{a^{2}}{2b}\bigl(c_{1}+a_{1}x+b_{1}y\bigr)^{2}
+ a\bigl(c_{1}+a_{1}x+b_{1}y\bigr)
   \bigl(x\alpha_{1}+y\beta_{1}+\gamma_{1}\bigr)
 \nonumber\\
&\quad + \frac{b}{2}\bigl(x\alpha_{1}+y\beta_{1}+\gamma_{1}\bigr)^{2}
- \frac{3}{2}\bigl(c_{1}+a_{1}x+b_{1}y\bigr)^{2}
  \bigl(x\alpha_{1}+y\beta_{1}+\gamma_{1}\bigr)^{2} \nonumber\\
&\quad
+ \frac{1}{4}\bigl(x\alpha_{1}+y\beta_{1}+\gamma_{1}\bigr)^{4}.
\end{align}

The differential system ($\mathtt{N}_4^{1}$) becomes

\begin{align}
\dot{x}
&= \frac{1}{-b_{1}\alpha_{1}+a_{1}\beta_{1}}
\Bigl(
- b_{1} c c_{1}
- 3 c_{1}^{2}\beta_{1}\gamma_{1}
- 3 b_{1} c_{1}\gamma_{1}^{2}
- \beta_{1}\gamma_{1}^{3} 
+ \bigl(-3 a_{1} b_{1}\alpha_{1}^{2} \nonumber\\
&\quad
- 3 a_{1}^{2}\alpha_{1}\beta_{1}
- \alpha_{1}^{3}\beta_{1}\bigr)x^{3}
- 3\bigl(b_{1}^{2}\alpha_{1}^{2}
+ 4 a_{1} b_{1}\alpha_{1}\beta_{1}
+ a_{1}^{2}\beta_{1}^{2}
+ \alpha_{1}^{2}\beta_{1}^{2}\bigr)x^{2}y \nonumber\\
&\quad
- 3\bigl(3 b_{1}^{2}\alpha_{1}\beta_{1}
+ 3 a_{1} b_{1}\beta_{1}^{2}
+ \alpha_{1}\beta_{1}^{3}\bigr)xy^{2}
+ \bigl(-6 b_{1}^{2}\beta_{1}^{2}
- \beta_{1}^{4}\bigr)y^{3} \nonumber\\
&\quad
- 3\bigl(b_{1} c_{1}\alpha_{1}^{2}
+ 2 a_{1} c_{1}\alpha_{1}\beta_{1}
+ 2 a_{1} b_{1}\alpha_{1}\gamma_{1}
+ a_{1}^{2}\beta_{1}\gamma_{1}
+ \alpha_{1}^{2}\beta_{1}\gamma_{1}\bigr)x^{2} \nonumber\\
&\quad
- 6\bigl(2 b_{1} c_{1}\alpha_{1}\beta_{1}
+ a_{1} c_{1}\beta_{1}^{2}
+ b_{1}^{2}\alpha_{1}\gamma_{1}
+ 2 a_{1} b_{1}\beta_{1}\gamma_{1}
+ \alpha_{1}\beta_{1}^{2}\gamma_{1}\bigr)xy \nonumber\\
&\quad
- 3\bigl(3 b_{1} c_{1}\beta_{1}^{2}
+ 3 b_{1}^{2}\beta_{1}\gamma_{1}
+ \beta_{1}^{3}\gamma_{1}\bigr)y^{2}
+ \bigl(-a_{1} b_{1} c
- 3 c_{1}^{2}\alpha_{1}\beta_{1}
- 6 b_{1} c_{1}\alpha_{1}\gamma_{1} \nonumber\\
&\quad
- 6 a_{1} c_{1}\beta_{1}\gamma_{1}
- 3 a_{1} b_{1}\gamma_{1}^{2}
- 3 \alpha_{1}\beta_{1}\gamma_{1}^{2}\bigr)x 
+ \bigl(-b_{1}^{2} c
- 3 c_{1}^{2}\beta_{1}^{2}
- 12 b_{1} c_{1}\beta_{1}\gamma_{1}\nonumber\\
&\quad
- 3 b_{1}^{2}\gamma_{1}^{2}
- 3 \beta_{1}^{2}\gamma_{1}^{2}\bigr)y
\Bigr),
\end{align}

\begin{align}
\dot{y}&=\frac{1}{-b_1\alpha_1 + a_1\beta_1}\big(
 a_1 c c_1
 + 3 c_1^2 \alpha_1 \gamma_1
 + 3 a_1 c_1 \gamma_1^2
 + \alpha_1 \gamma_1^3  + \left(6 a_1^2 \alpha_1^2
 + \alpha_1^4\right)x^3 \nonumber\\
& + 3\left(3 a_1 b_1 \alpha_1^2
 + 3 a_1^2 \alpha_1 \beta_1
 + \alpha_1^3 \beta_1\right)x^2 y  + 3(b_1^2 \alpha_1^2
 + 4 a_1 b_1 \alpha_1 \beta_1
 + a_1^2 \beta_1^2 \nonumber\\
&
 + \alpha_1^2 \beta_1^2)x y^2 + \beta_1\left(3 b_1^2 \alpha_1
 + 3 a_1 b_1 \beta_1
 + \alpha_1 \beta_1^2\right)y^3 + 3(3 a_1 c_1 \alpha_1^2
 + 3 a_1^2 \alpha_1 \gamma_1\nonumber\\
& 
 + \alpha_1^3 \gamma_1)x^2 + 6\left(b_1 c_1 \alpha_1^2
 + 2 a_1 c_1 \alpha_1 \beta_1
 + 2 a_1 b_1 \alpha_1 \gamma_1
 + a_1^2 \beta_1 \gamma_1
 + \alpha_1^2 \beta_1 \gamma_1\right)x y \nonumber\\
& + 3\left(2 b_1 c_1 \alpha_1 \beta_1
 + a_1 c_1 \beta_1^2
 + b_1^2 \alpha_1 \gamma_1
 + 2 a_1 b_1 \beta_1 \gamma_1
 + \alpha_1 \beta_1^2 \gamma_1\right)y^2 \nonumber\\
& + \left(a_1^2 c
 + 3 c_1^2 \alpha_1^2
 + 12 a_1 c_1 \alpha_1 \gamma_1
 + 3 a_1^2 \gamma_1^2
 + 3 \alpha_1^2 \gamma_1^2\right)x \nonumber\\
& + \left(a_1 b_1 c
 + 3 c_1^2 \alpha_1 \beta_1
 + 6 b_1 c_1 \alpha_1 \gamma_1
 + 6 a_1 c_1 \beta_1 \gamma_1
 + 3 a_1 b_1 \gamma_1^2
 + 3 \alpha_1 \beta_1 \gamma_1^2\right)y
\big),
\end{align}

with the first integral
\begin{align}
H_4^{1}(x,y)=&-\frac{1}{2}c\left(c_{1}+a_{1}x+b_{1}y\right)^{2}
-\frac{3}{2}\left(c_{1}+a_{1}x+b_{1}y\right)^{2}
\left(x\alpha_{1}+y\beta_{1}+\gamma_{1}\right)^{2} \nonumber\\
&\quad
-\frac{1}{4}\left(x\alpha_{1}+y\beta_{1}+\gamma_{1}\right)^{4}.
\end{align}

The differential system ($\mathtt{N}_4^{2}$) becomes

\begin{equation}
\begin{aligned}
\dot{x}
&= \frac{1}{b\left(b_{1}\alpha_{1}-a_{1}\beta_{1}\right)}
\Bigl(
- a^{2} b_{1} c_{1}
- a b c_{1}\beta_{1}
- a b b_{1}\gamma_{1}
- b^{2}\beta_{1}\gamma_{1}
+ 3 b c_{1}^{2}\beta_{1}\gamma_{1} \\
&\quad
+ 3 b b_{1} c_{1}\gamma_{1}^{2}
+ b\beta_{1}\gamma_{1}^{3}
+ b\alpha_{1}\bigl(3 a_{1} b_{1}\alpha_{1}
+ 3 a_{1}^{2}\beta_{1}
+ \alpha_{1}^{2}\beta_{1}\bigr)x^{3}
+ 3\bigl(b b_{1}^{2}\alpha_{1}^{2} \\
&\quad
+ 4 a_{1} b b_{1}\alpha_{1}\beta_{1}
+ a_{1}^{2} b\beta_{1}^{2}
+ b\alpha_{1}^{2}\beta_{1}^{2}\bigr)x^{2}y
+ 3\bigl(3 b b_{1}^{2}\alpha_{1}\beta_{1}
+ 3 a_{1} b b_{1}\beta_{1}^{2} \\
&\quad
+ b\alpha_{1}\beta_{1}^{3}\bigr)xy^{2}
+ b\beta_{1}^{2}\bigl(6 b_{1}^{2}
+ \beta_{1}^{2}\bigr)y^{3}
+ 3\bigl(b b_{1} c_{1}\alpha_{1}^{2}
+ 2 a_{1} b c_{1}\alpha_{1}\beta_{1}
+ 2 a_{1} b b_{1}\alpha_{1}\gamma_{1} \\
&\quad
+ a_{1}^{2} b\beta_{1}\gamma_{1}
+ b\alpha_{1}^{2}\beta_{1}\gamma_{1}\bigr)x^{2} 
+ 6\bigl(2 b b_{1} c_{1}\alpha_{1}\beta_{1}
+ a_{1} b c_{1}\beta_{1}^{2}
+ b b_{1}^{2}\alpha_{1}\gamma_{1}\\
&\quad
+ 2 a_{1} b b_{1}\beta_{1}\gamma_{1}
+ b\alpha_{1}\beta_{1}^{2}\gamma_{1}\bigr)xy 
+ 3\bigl(3 b b_{1} c_{1}\beta_{1}^{2}
+ 3 b b_{1}^{2}\beta_{1}\gamma_{1}
+ b\beta_{1}^{3}\gamma_{1}\bigr)y^{2} \\
&\quad
+ \bigl(-a^{2} a_{1} b_{1}
- a b b_{1}\alpha_{1}
- a a_{1} b\beta_{1}
- b^{2}\alpha_{1}\beta_{1}
+ 3 b c_{1}^{2}\alpha_{1}\beta_{1}
+ 6 b b_{1} c_{1}\alpha_{1}\gamma_{1} \\
&\quad
+ 6 a_{1} b c_{1}\beta_{1}\gamma_{1}
+ 3 a_{1} b b_{1}\gamma_{1}^{2}
+ 3 b\alpha_{1}\beta_{1}\gamma_{1}^{2}\bigr)x
+ \bigl(-a^{2} b_{1}^{2}
- 2 a b b_{1}\beta_{1}
- b^{2}\beta_{1}^{2}\\
&\quad
+ 3 b c_{1}^{2}\beta_{1}^{2}
+ 12 b b_{1} c_{1}\beta_{1}\gamma_{1}
+ 3 b b_{1}^{2}\gamma_{1}^{2}
+ 3 b\beta_{1}^{2}\gamma_{1}^{2}\bigr)y
\Bigr),
\end{aligned}
\end{equation}

\begin{equation}
\begin{aligned}
\dot{y}
&= \frac{1}{b\left(b_{1}\alpha_{1}-a_{1}\beta_{1}\right)}
\big(
a^{2} a_{1} c_{1}
+ a b c_{1}\alpha_{1}
+ a a_{1} b\gamma_{1}
+ b^{2}\alpha_{1}\gamma_{1}
- 3 b c_{1}^{2}\alpha_{1}\gamma_{1}\\
&\quad
- 3 a_{1} b c_{1}\gamma_{1}^{2}
- b\alpha_{1}\gamma_{1}^{3} 
+ \bigl(-6 a_{1}^{2} b\alpha_{1}^{2}
- b\alpha_{1}^{4}\bigr)x^{3}
- 3\bigl(3 a_{1} b b_{1}\alpha_{1}^{2}
+ 3 a_{1}^{2} b\alpha_{1}\beta_{1}\\
&\quad
+ b\alpha_{1}^{3}\beta_{1}\bigr)x^{2}y 
- 3 b\bigl(b_{1}^{2}\alpha_{1}^{2}
+ 4 a_{1} b_{1}\alpha_{1}\beta_{1}
+ a_{1}^{2}\beta_{1}^{2}
+ \alpha_{1}^{2}\beta_{1}^{2}\bigr)xy^{2} \\
&\quad
+ \bigl(-3 b b_{1}^{2}\alpha_{1}\beta_{1}
- 3 a_{1} b b_{1}\beta_{1}^{2}
- b\alpha_{1}\beta_{1}^{3}\bigr)y^{3}
- 3\bigl(3 a_{1} b c_{1}\alpha_{1}^{2}
+ 3 a_{1}^{2} b\alpha_{1}\gamma_{1}  \\
&\quad
+ b\alpha_{1}^{3}\gamma_{1}\bigr)x^{2}
- 6 b\bigl(b_{1} c_{1}\alpha_{1}^{2}
+ 2 a_{1} c_{1}\alpha_{1}\beta_{1}
+ 2 a_{1} b_{1}\alpha_{1}\gamma_{1}
+ a_{1}^{2}\beta_{1}\gamma_{1}
+ \alpha_{1}^{2}\beta_{1}\gamma_{1}\bigr)xy \\
&\quad
- 3 b\bigl(2 b_{1} c_{1}\alpha_{1}\beta_{1}
+ a_{1} c_{1}\beta_{1}^{2}
+ b_{1}^{2}\alpha_{1}\gamma_{1}
+ 2 a_{1} b_{1}\beta_{1}\gamma_{1}
+ \alpha_{1}\beta_{1}^{2}\gamma_{1}\bigr)y^{2} \\
&\quad
+ \bigl(a^{2} a_{1}^{2}
+ 2 a a_{1} b\alpha_{1}
+ b^{2}\alpha_{1}^{2}
- 3 b c_{1}^{2}\alpha_{1}^{2}
- 12 a_{1} b c_{1}\alpha_{1}\gamma_{1}
- 3 a_{1}^{2} b\gamma_{1}^{2}
- 3 b\alpha_{1}^{2}\gamma_{1}^{2}\bigr)x \\
&\quad
+ \bigl(a^{2} a_{1} b_{1}
+ a b b_{1}\alpha_{1}
+ a a_{1} b\beta_{1}
+ b^{2}\alpha_{1}\beta_{1}
- 3 b c_{1}^{2}\alpha_{1}\beta_{1}
- 6 b b_{1} c_{1}\alpha_{1}\gamma_{1}\\
&\quad
- 6 a_{1} b c_{1}\beta_{1}\gamma_{1}
- 3 a_{1} b b_{1}\gamma_{1}^{2}
- 3 b\alpha_{1}\beta_{1}\gamma_{1}^{2}\bigr)y
\big),
\end{aligned}
\end{equation}

with the first integral
\begin{equation}
\begin{aligned}
H_4^{2}(x,y)
&= \frac{a^{2}}{2b}\left(c_{1}+a_{1}x+b_{1}y\right)^{2}
+ a\left(c_{1}+a_{1}x+b_{1}y\right)
\left(x\alpha_{1}+y\beta_{1}+\gamma_{1}\right)\\
&\quad
+ \frac{b}{2}\left(x\alpha_{1}+y\beta_{1}+\gamma_{1}\right)^{2} 
- \frac{3}{2}\left(c_{1}+a_{1}x+b_{1}y\right)^{2}
\left(x\alpha_{1}+y\beta_{1}+\gamma_{1}\right)^{2}
\\
&\quad
- \frac{1}{4}\left(x\alpha_{1}+y\beta_{1}+\gamma_{1}\right)^{4}.
\end{aligned}
\end{equation}


The differential system ($\mathtt{N}_5^{1}$) becomes

\begin{equation}
\begin{aligned}
\dot{x}
&=\frac{1}{-b_1\alpha_1+a_1\beta_1}\big(
 - b_1 c c_1
 - b_1 c_1^3
 + \beta_1 \gamma_1^3
 - 3 c_1^2 \beta_1 \gamma_1 \mu
 - 3 b_1 c_1 \gamma_1^2 \mu  
 + ( - a_1^3 b_1 \\
&\quad
 + \alpha_1^3 \beta_1
 - 3 a_1 b_1 \alpha_1^2 \mu
 - 3 a_1^2 \alpha_1 \beta_1 \mu
 )x^3 
 - 3(
 a_1^2 b_1^2
 - \alpha_1^2 \beta_1^2
 + b_1^2 \alpha_1^2 \mu  \\
&\quad
 + 4 a_1 b_1 \alpha_1 \beta_1 \mu
 + a_1^2 \beta_1^2 \mu
 )x^2 y  
 - 3\left(
 a_1 b_1^3
 - \alpha_1 \beta_1^3
 + 3 b_1^2 \alpha_1 \beta_1 \mu
 + 3 a_1 b_1 \beta_1^2 \mu
 \right)x y^2  \\
&\quad
 + \left(
 - b_1^4
 + \beta_1^4
 - 6 b_1^2 \beta_1^2 \mu
 \right)y^3  
 - 3(
 a_1^2 b_1 c_1
 - \alpha_1^2 \beta_1 \gamma_1
 + b_1 c_1 \alpha_1^2 \mu \\
&\quad
 + 2 a_1 c_1 \alpha_1 \beta_1 \mu
 + 2 a_1 b_1 \alpha_1 \gamma_1 \mu
 + a_1^2 \beta_1 \gamma_1 \mu
)x^2   - 6( a_1 b_1^2 c_1
 - \alpha_1 \beta_1^2 \gamma_1\\
&\quad
 + 2 b_1 c_1 \alpha_1 \beta_1 \mu
 + a_1 c_1 \beta_1^2 \mu
 + b_1^2 \alpha_1 \gamma_1 \mu
 + 2 a_1 b_1 \beta_1 \gamma_1 \mu
)x y  
 - 3(
 b_1^3 c_1
 - \beta_1^3 \gamma_1 \\
&\quad
 + 3 b_1 c_1 \beta_1^2 \mu
 + 3 b_1^2 \beta_1 \gamma_1 \mu
)y^2 
 + (
 - a_1 b_1 c
 - 3 a_1 b_1 c_1^2
 + 3 \alpha_1 \beta_1 \gamma_1^2
 - 3 c_1^2 \alpha_1 \beta_1 \mu \\
&\quad
 - 6 b_1 c_1 \alpha_1 \gamma_1 \mu
 - 6 a_1 c_1 \beta_1 \gamma_1 \mu
 - 3 a_1 b_1 \gamma_1^2 \mu
)x  
 + (
 - b_1^2 c
 - 3 b_1^2 c_1^2
 + 3 \beta_1^2 \gamma_1^2\\
&\quad
 - 3 c_1^2 \beta_1^2 \mu
 - 12 b_1 c_1 \beta_1 \gamma_1 \mu
 - 3 b_1^2 \gamma_1^2 \mu
 )y
\big),
\end{aligned}
\end{equation}


\begin{equation}
\begin{aligned}
\dot{y}
&=\frac{1}{-b_1\alpha_1+a_1\beta_1}\big(
 a_1 c c_1
 + a_1 c_1^3
 - \alpha_1 \gamma_1^3
 + 3 c_1^2 \alpha_1 \gamma_1 \mu
 + 3 a_1 c_1 \gamma_1^2 \mu  
 + (
 a_1^4
 - \alpha_1^4 \\
&\quad
 + 6 a_1^2 \alpha_1^2 \mu
)x^3 + 3\left(
 a_1^3 b_1
 - \alpha_1^3 \beta_1
 + 3 a_1 b_1 \alpha_1^2 \mu
 + 3 a_1^2 \alpha_1 \beta_1 \mu
 \right)x^2 y 
 + 3(
 a_1^2 b_1^2 \\
&\quad
 - \alpha_1^2 \beta_1^2
 + b_1^2 \alpha_1^2 \mu
 + 4 a_1 b_1 \alpha_1 \beta_1 \mu
 + a_1^2 \beta_1^2 \mu
)x y^2  
 + (
 a_1 b_1^3
 - \alpha_1 \beta_1^3
 + 3 b_1^2 \alpha_1 \beta_1 \mu \\
&\quad
 + 3 a_1 b_1 \beta_1^2 \mu
 )y^3 
 + 3\left(
 a_1^3 c_1
 - \alpha_1^3 \gamma_1
 + 3 a_1 c_1 \alpha_1^2 \mu
 + 3 a_1^2 \alpha_1 \gamma_1 \mu
 \right)x^2  \\
&\quad
 + 6(
 a_1^2 b_1 c_1
 - \alpha_1^2 \beta_1 \gamma_1
 + b_1 c_1 \alpha_1^2 \mu
 + 2 a_1 c_1 \alpha_1 \beta_1 \mu
 + 2 a_1 b_1 \alpha_1 \gamma_1 \mu \\
&\quad
 + a_1^2 \beta_1 \gamma_1 \mu
)x y 
 + 3(
 a_1 b_1^2 c_1
 - \alpha_1 \beta_1^2 \gamma_1
 + 2 b_1 c_1 \alpha_1 \beta_1 \mu
 + a_1 c_1 \beta_1^2 \mu
 + b_1^2 \alpha_1 \gamma_1 \mu \\
&\quad
 + 2 a_1 b_1 \beta_1 \gamma_1 \mu
)y^2 
 + (
 a_1^2 c
 + 3 a_1^2 c_1^2
 - 3 \alpha_1^2 \gamma_1^2
 + 3 c_1^2 \alpha_1^2 \mu
 + 12 a_1 c_1 \alpha_1 \gamma_1 \mu \\
&\quad
 + 3 a_1^2 \gamma_1^2 \mu
 )x 
 + (
 a_1 b_1 c
 + 3 a_1 b_1 c_1^2
 - 3 \alpha_1 \beta_1 \gamma_1^2
 + 3 c_1^2 \alpha_1 \beta_1 \mu
 + 6 b_1 c_1 \alpha_1 \gamma_1 \mu
 \\
&\quad
 + 6 a_1 c_1 \beta_1 \gamma_1 \mu
 + 3 a_1 b_1 \gamma_1^2 \mu
 )y
\big),
\end{aligned}
\end{equation}

with the first integral
\begin{align}
H_5^{1}(x,y)&=-\frac{1}{2}c\left(c_{1}+a_{1}x+b_{1}y\right)^{2}
+\frac{1}{4}\big(
-\left(c_{1}+a_{1}x+b_{1}y\right)^{4}
+(x\alpha_{1}+y\beta_{1}\nonumber\\
&\quad +\gamma_{1})^{4}
\big) 
-\frac{3}{2}\mu\left(c_{1}+a_{1}x+b_{1}y\right)^{2}
\left(x\alpha_{1}+y\beta_{1}+\gamma_{1}\right)^{2}.
\end{align}

The differential system ($\mathtt{N}_5^{2}$) becomes
\begin{equation}
\begin{aligned}
\dot{x}
=&\frac{1}{b(b_1\alpha_1-a_1\beta_1)}\big(
-a^{2}b_{1}c_{1}+bb_{1}c_{1}^{3}-abc_{1}\beta_{1}-abb_{1}\gamma_{1}
-b^{2}\beta_{1}\gamma_{1}-b\beta_{1}\gamma_{1}^{3} \\
&\quad
+3bc_{1}^{2}\beta_{1}\gamma_{1}\mu
+3bb_{1}c_{1}\gamma_{1}^{2}\mu
+b\!\left(a_{1}^{3}b_{1}-\alpha_{1}^{3}\beta_{1}
+3a_{1}b_{1}\alpha_{1}^{2}\mu
+3a_{1}^{2}\alpha_{1}\beta_{1}\mu\right)x^{3} \\
&\quad
+3\!\left(a_{1}^{2}bb_{1}^{2}-b\alpha_{1}^{2}\beta_{1}^{2}
+bb_{1}^{2}\alpha_{1}^{2}\mu
+4a_{1}bb_{1}\alpha_{1}\beta_{1}\mu
+a_{1}^{2}b\beta_{1}^{2}\mu\right)x^{2}y \\
&\quad
+3\!\left(a_{1}bb_{1}^{3}-b\alpha_{1}\beta_{1}^{3}
+3bb_{1}^{2}\alpha_{1}\beta_{1}\mu
+3a_{1}bb_{1}\beta_{1}^{2}\mu\right)xy^{2}
+b(b_{1}^{4}-\beta_{1}^{4} \\
&\quad
+6b_{1}^{2}\beta_{1}^{2}\mu)y^{3} 
+3(a_{1}^{2}bb_{1}c_{1}-b\alpha_{1}^{2}\beta_{1}\gamma_{1}
+bb_{1}c_{1}\alpha_{1}^{2}\mu
+2a_{1}bc_{1}\alpha_{1}\beta_{1}\mu\\
&\quad
+2a_{1}bb_{1}\alpha_{1}\gamma_{1}\mu
+a_{1}^{2}b\beta_{1}\gamma_{1}\mu)x^{2} 
+6(a_{1}bb_{1}^{2}c_{1}-b\alpha_{1}\beta_{1}^{2}\gamma_{1} \\
&\quad
+2bb_{1}c_{1}\alpha_{1}\beta_{1}\mu
+a_{1}bc_{1}\beta_{1}^{2}\mu
+bb_{1}^{2}\alpha_{1}\gamma_{1}\mu
+2a_{1}bb_{1}\beta_{1}\gamma_{1}\mu)xy
+3(bb_{1}^{3}c_{1}  \\
&\quad
-b\beta_{1}^{3}\gamma_{1}
+3bb_{1}c_{1}\beta_{1}^{2}\mu
+3bb_{1}^{2}\beta_{1}\gamma_{1}\mu)y^{2}
+(-a^{2}a_{1}b_{1}+3a_{1}bb_{1}c_{1}^{2}
-abb_{1}\alpha_{1}\\
&\quad
-aa_{1}b\beta_{1}
-b^{2}\alpha_{1}\beta_{1}-3b\alpha_{1}\beta_{1}\gamma_{1}^{2} 
+3bc_{1}^{2}\alpha_{1}\beta_{1}\mu
+6bb_{1}c_{1}\alpha_{1}\gamma_{1}\mu \\
&\quad
+6a_{1}bc_{1}\beta_{1}\gamma_{1}\mu
+3a_{1}bb_{1}\gamma_{1}^{2}\mu)x
+(-a^{2}b_{1}^{2}+3bb_{1}^{2}c_{1}^{2}
-2abb_{1}\beta_{1}-b^{2}\beta_{1}^{2} \\
&\quad
-3b\beta_{1}^{2}\gamma_{1}^{2}
+3bc_{1}^{2}\beta_{1}^{2}\mu
+12bb_{1}c_{1}\beta_{1}\gamma_{1}\mu
+3bb_{1}^{2}\gamma_{1}^{2}\mu)y
\big),
\end{aligned}
\end{equation}

\begin{equation}
\begin{aligned}
\dot{y}
=&\frac{1}{b(b_1\alpha_1-a_1\beta_1)}\big(
a^2a_1c_1-a_1bc_1^3+abc_1\alpha_1+aa_1b\gamma_1
+b^2\alpha_1\gamma_1+b\alpha_1\gamma_1^3 \\
&\quad
-3bc_1^2\alpha_1\gamma_1\mu
-3a_1bc_1\gamma_1^2\mu
+\big(-a_1^4b+b\alpha_1^4-6a_1^2b\alpha_1^2\mu\big)x^3
-3\big(a_1^3bb_1  \\
&\quad
-b\alpha_1^3\beta_1
+3a_1bb_1\alpha_1^2\mu
+3a_1^2b\alpha_1\beta_1\mu\big)x^2y 
-3b\big(a_1^2b_1^2-\alpha_1^2\beta_1^2 
+b_1^2\alpha_1^2\mu\\
&\quad
+4a_1b_1\alpha_1\beta_1\mu
+a_1^2\beta_1^2\mu\big)xy^2
+\big(-a_1bb_1^3+b\alpha_1\beta_1^3
-3bb_1^2\alpha_1\beta_1\mu\\
&\quad
-3a_1bb_1\beta_1^2\mu\big)y^3 
-3\big(a_1^3bc_1-b\alpha_1^3\gamma_1
+3a_1bc_1\alpha_1^2\mu
+3a_1^2b\alpha_1\gamma_1\mu\big)x^2 \\
&\quad
-6b\big(a_1^2b_1c_1-\alpha_1^2\beta_1\gamma_1
+b_1c_1\alpha_1^2\mu
+2a_1c_1\alpha_1\beta_1\mu
+2a_1b_1\alpha_1\gamma_1\mu \\
&\quad
+a_1^2\beta_1\gamma_1\mu\big)xy 
-3b\big(a_1b_1^2c_1-\alpha_1\beta_1^2\gamma_1
+2b_1c_1\alpha_1\beta_1\mu
+a_1c_1\beta_1^2\mu \\
&\quad
+b_1^2\alpha_1\gamma_1\mu
+2a_1b_1\beta_1\gamma_1\mu\big)y^2
+\big(a^2a_1^2-3a_1^2bc_1^2
+2aa_1b\alpha_1+b^2\alpha_1^2 \\
&\quad
+3b\alpha_1^2\gamma_1^2
-3bc_1^2\alpha_1^2\mu 
-12a_1bc_1\alpha_1\gamma_1\mu
-3a_1^2b\gamma_1^2\mu\big)x 
+\big(a^2a_1b_1 \\
&\quad
-3a_1bb_1c_1^2
+abb_1\alpha_1+aa_1b\beta_1
+b^2\alpha_1\beta_1
+3b\alpha_1\beta_1\gamma_1^2  \\
&\quad
-3bc_1^2\alpha_1\beta_1\mu
-6bb_1c_1\alpha_1\gamma_1\mu
-6a_1bc_1\beta_1\gamma_1\mu
-3a_1bb_1\gamma_1^2\mu\big)y
\big),
\end{aligned}
\end{equation}

with the first integral
\begin{equation}
\begin{aligned}
H_5^{2}(x,y)=& \frac{a^{2}}{2b}\left(c_{1}+a_{1}x+b_{1}y\right)^{2}
+a\left(c_{1}+a_{1}x+b_{1}y\right)
\left(x\alpha_{1}+y\beta_{1}+\gamma_{1}\right)\\
&
+\frac{b}{2}\left(x\alpha_{1}+y\beta_{1}+\gamma_{1}\right)^{2} 
+\frac{1}{4}\big(
-\left(c_{1}+a_{1}x+b_{1}y\right)^{4}
+(x\alpha_{1}+y\beta_{1}\\
&+\gamma_{1})^{4}
\big)
-\frac{3}{2}\mu\left(c_{1}+a_{1}x+b_{1}y\right)^{2}
\left(x\alpha_{1}+y\beta_{1}+\gamma_{1}\right)^{2}.
\end{aligned}
\end{equation}

The differential system ($\mathtt{N}_6^{1}$) becomes
\begin{equation}
\begin{aligned}
\dot{x}=&\frac{1}{-b_1\alpha_1+a_1\beta_1}\big(- b_1 c c_1 - b_1 c_1^3 - \beta_1 \gamma_1^3 - 3 c_1^2 \beta_1 \gamma_1 \mu - 3 b_1 c_1 \gamma_1^2 \mu+ (-a_1^3 b_1\nonumber\\
& - \alpha_1^3 \beta_1 - 3 a_1 b_1 \alpha_1^2 \mu 
- 3 a_1^2 \alpha_1 \beta_1 \mu)\,x^3 - 3\big(a_1^2 b_1^2 + \alpha_1^2 \beta_1^2
+ b_1^2 \alpha_1^2 \mu \nonumber\\
&+ 4 a_1 b_1 \alpha_1 \beta_1 \mu
+ a_1^2 \beta_1^2 \mu\big)x^2 y- 3\big(a_1 b_1^3 + \alpha_1 \beta_1^3
+ 3 b_1^2 \alpha_1 \beta_1 \mu \nonumber\\
& + 3 a_1 b_1 \beta_1^2 \mu\big)x y^2 
+ (-b_1^4 - \beta_1^4 - 6 b_1^2 \beta_1^2 \mu)\,y^3
- 3\big(a_1^2 b_1 c_1 + \alpha_1^2 \beta_1 \gamma_1  \nonumber\\
&
+ b_1 c_1 \alpha_1^2 \mu 
+ 2 a_1 c_1 \alpha_1 \beta_1 \mu + 2 a_1 b_1 \alpha_1 \gamma_1 \mu
+ a_1^2 \beta_1 \gamma_1 \mu\big)x^2 - 6\big(a_1 b_1^2 c_1  \nonumber\\
& 
+ \alpha_1 \beta_1^2 \gamma_1 + 2 b_1 c_1 \alpha_1 \beta_1 \mu
+ a_1 c_1 \beta_1^2 \mu+ b_1^2 \alpha_1 \gamma_1 \mu
+ 2 a_1 b_1 \beta_1 \gamma_1 \mu\big)x y \nonumber\\
&
- 3\big(b_1^3 c_1 + \beta_1^3 \gamma_1
+ 3 b_1 c_1 \beta_1^2 \mu
+ 3 b_1^2 \beta_1 \gamma_1 \mu\big)y^2+ 
\big(-a_1 b_1 c - 3 a_1 b_1 c_1^2  \nonumber\\
&
- 3 \alpha_1 \beta_1 \gamma_1^2
- 3 c_1^2 \alpha_1 \beta_1 \mu - 6 b_1 c_1 \alpha_1 \gamma_1 \mu
- 6 a_1 c_1 \beta_1 \gamma_1 \mu
- 3 a_1 b_1 \gamma_1^2 \mu\big)x \nonumber\\
&+ \big(-b_1^2 c - 3 b_1^2 c_1^2
- 3 \beta_1^2 \gamma_1^2
- 3 c_1^2 \beta_1^2 \mu  - 12 b_1 c_1 \beta_1 \gamma_1 \mu
- 3 b_1^2 \gamma_1^2 \mu\big)y
\big),
\end{aligned}
\end{equation}

\begin{equation}
\begin{aligned}
\dot{y}&=\frac{1}{-b_1\alpha_1+a_1\beta_1}\big( a_1 c c_1 + a_1 c_1^3 + \alpha_1 \gamma_1^3 + 3 c_1^2 \alpha_1 \gamma_1 \mu + 3 a_1 c_1 \gamma_1^2 \mu + (a_1^4 + \alpha_1^4 \\
& \quad 
+ 6 a_1^2 \alpha_1^2 \mu)\,x^3 + 3\big(a_1^3 b_1 + \alpha_1^3 \beta_1 + 3 a_1 b_1 \alpha_1^2 \mu + 3 a_1^2 \alpha_1 \beta_1 \mu\big)x^2 y + 3\big(a_1^2 b_1^2
\\
& \quad 
+ \alpha_1^2 \beta_1^2
+ b_1^2 \alpha_1^2 \mu + 4 a_1 b_1 \alpha_1 \beta_1 \mu
+ a_1^2 \beta_1^2 \mu\big)x y^2  + \big(a_1 b_1^3 + \alpha_1 \beta_1^3  + 3 b_1^2 \alpha_1 \beta_1 \mu 
\\
& \quad 
+ 3 a_1 b_1 \beta_1^2 \mu\big)y^3 + 3\big(a_1^3 c_1 + \alpha_1^3 \gamma_1
+ 3 a_1 c_1 \alpha_1^2 \mu + 3 a_1^2 \alpha_1 \gamma_1 \mu\big)x^2  + 6\big(a_1^2 b_1 c_1 \\
& \quad 
+ \alpha_1^2 \beta_1 \gamma_1
+ b_1 c_1 \alpha_1^2 \mu + 2 a_1 c_1 \alpha_1 \beta_1 \mu  + 2 a_1 b_1 \alpha_1 \gamma_1 \mu + a_1^2 \beta_1 \gamma_1 \mu\big)x y \\
& \quad  + 3\big(a_1 b_1^2 c_1 + \alpha_1 \beta_1^2 \gamma_1
+ 2 b_1 c_1 \alpha_1 \beta_1 \mu + a_1 c_1 \beta_1^2 \mu  + b_1^2 \alpha_1 \gamma_1 \mu + 2 a_1 b_1 \beta_1 \gamma_1 \mu\big)y^2 \\
& \quad + \big(a_1^2 c + 3 a_1^2 c_1^2 + 3 \alpha_1^2 \gamma_1^2
+ 3 c_1^2 \alpha_1^2 \mu + 12 a_1 c_1 \alpha_1 \gamma_1 \mu
+ 3 a_1^2 \gamma_1^2 \mu\big)x \\
& \quad  + \big(a_1 b_1 c + 3 a_1 b_1 c_1^2
+ 3 \alpha_1 \beta_1 \gamma_1^2
+ 3 c_1^2 \alpha_1 \beta_1 \mu  + 6 b_1 c_1 \alpha_1 \gamma_1 \mu
+ 6 a_1 c_1 \beta_1 \gamma_1 \mu \\
& \quad 
+ 3 a_1 b_1 \gamma_1^2 \mu\big)y
\Big),
\end{aligned}
\end{equation}
with the first integral
\begin{equation}
\begin{aligned}
H_6^{1}(x,y)=&-\frac{1}{2}c\left(c_{1}+a_{1}x+b_{1}y\right)^{2}
+\frac{1}{4}\big(
-\left(c_{1}+a_{1}x+b_{1}y\right)^{4}
-(x\alpha_{1}+y\beta_{1} \nonumber\\
&\quad +\gamma_{1})^{4}
\big) 
-\frac{3}{2}\mu\left(c_{1}+a_{1}x+b_{1}y\right)^{2}
\left(x\alpha_{1}+y\beta_{1}+\gamma_{1}\right)^{2}.
\end{aligned}
\end{equation}
The differential system ($\mathtt{N}_6^{2}$) becomes
\begin{equation}
\begin{aligned}
\dot{x}&=\frac{1}{b(b_1\alpha_1-a_1\beta_1)}\big(
-a^2 b_1 c_1 + b b_1 c_1^3
- a b c_1 \beta_1 - a b b_1 \gamma_1
- b^2 \beta_1 \gamma_1
+ b \beta_1 \gamma_1^3 \\ & \quad + 3 b c_1^2 \beta_1 \gamma_1 \mu
+ 3 b b_1 c_1 \gamma_1^2 \mu + b\big(a_1^3 b_1 + \alpha_1^3 \beta_1
+ 3 a_1 b_1 \alpha_1^2 \mu
+ 3 a_1^2 \alpha_1 \beta_1 \mu\big)x^3 \\ & \quad+ 3\big(a_1^2 b b_1^2 + b \alpha_1^2 \beta_1^2 + b b_1^2 \alpha_1^2 \mu
+ 4 a_1 b b_1 \alpha_1 \beta_1 \mu + a_1^2 b \beta_1^2 \mu\big)x^2 y \\ & \quad+ 3\big(a_1 b b_1^3 + b \alpha_1 \beta_1^3
+ 3 b b_1^2 \alpha_1 \beta_1 \mu
+ 3 a_1 b b_1 \beta_1^2 \mu\big)x y^2+ b\big(b_1^4 + \beta_1^4 + 6 b_1^2 \beta_1^2 \mu\big)y^3  
\\ & \quad+ 3\big(a_1^2 b b_1 c_1 + b \alpha_1^2 \beta_1 \gamma_1
+ b b_1 c_1 \alpha_1^2 \mu
+ 2 a_1 b c_1 \alpha_1 \beta_1  + 2 a_1 b b_1 \alpha_1 \gamma_1 \mu\\ & \quad
+ a_1^2 b \beta_1 \gamma_1 \mu\big)x^2 
+ 6\big(a_1 b b_1^2 c_1 + b \alpha_1 \beta_1^2 \gamma_1
+ 2 b b_1 c_1 \alpha_1 \beta_1 \mu
+ a_1 b c_1 \beta_1^2 \mu\\ & \quad + b b_1^2 \alpha_1 \gamma_1 \mu
+ 2 a_1 b b_1 \beta_1 \gamma_1 \mu\big)x y+ 3\big(b b_1^3 c_1 + b \beta_1^3 \gamma_1
+ 3 b b_1 c_1 \beta_1^2 \mu  \\ & \quad
+ 3 b b_1^2 \beta_1 \gamma_1 \mu\big)y^2 + \big(-a^2 a_1 b_1  + 3 a_1 b b_1 c_1^2 - a b b_1 \alpha_1- a a_1 b \beta_1
- b^2 \alpha_1 \beta_1 \\ & \quad + 3 b \alpha_1 \beta_1 \gamma_1^2
+ 3 b c_1^2 \alpha_1 \beta_1 \mu + 6 b b_1 c_1 \alpha_1 \gamma_1 \mu  + 6 a_1 b c_1 \beta_1 \gamma_1 \mu + 3 a_1 b b_1 \gamma_1^2 \mu\big)x \\ & \quad + \big(-a^2 b_1^2 + 3 b b_1^2 c_1^2
- 2 a b b_1 \beta_1 - b^2 \beta_1^2
+ 3 b \beta_1^2 \gamma_1^2  + 3 b c_1^2 \beta_1^2 \mu
+ 12 b b_1 c_1 \beta_1 \gamma_1 \mu \\ & \quad
+ 3 b b_1^2 \gamma_1^2 \mu\big)y
\big),
\end{aligned}
\end{equation}

\begin{equation}
\begin{aligned}
\dot{y}&=\frac{1}{b\,(b_1\alpha_1-a_1\beta_1)}
\big(  a^{2}a_{1}c_{1}-a_{1}b c_{1}^{3}
+ab c_{1}\alpha_{1}+a a_{1}b\gamma_{1}
+b^{2}\alpha_{1}\gamma_{1}
-b\alpha_{1}\gamma_{1}^{3}\\
& \quad -3b c_{1}^{2}\alpha_{1}\gamma_{1}\mu  
-3a_{1}b c_{1}\gamma_{1}^{2}\mu +(-a_{1}^{4}b-b\alpha_{1}^{4}
-6a_{1}^{2}b\alpha_{1}^{2}\mu)x^{3} -3(a_{1}^{3}b b_{1}
+b\alpha_{1}^{3}\beta_{1}
\\ & \quad  +3a_{1}b b_{1}\alpha_{1}^{2}\mu  
+3a_{1}^{2}b\alpha_{1}\beta_{1}\mu)x^{2}y-3b(a_{1}^{2}b_{1}^{2} +\alpha_{1}^{2}\beta_{1}^{2} +b_{1}^{2}\alpha_{1}^{2}\mu
+4a_{1}b_{1}\alpha_{1}\beta_{1}\mu 
\\
& \quad +a_{1}^{2}\beta_{1}^{2}\mu)xy^{2} +(-a_{1}b b_{1}^{3} -b\alpha_{1}\beta_{1}^{3}
-3b b_{1}^{2}\alpha_{1}\beta_{1}\mu -3a_{1}b b_{1}\beta_{1}^{2}\mu)y^{3} -3(a_{1}^{3}b c_{1}\\
& \quad
+b\alpha_{1}^{3}\gamma_{1}
+3a_{1}b c_{1}\alpha_{1}^{2}\mu 
+3a_{1}^{2}b\alpha_{1}\gamma_{1}\mu)x^{2}
-6b(a_{1}^{2}b_{1}c_{1}
+\alpha_{1}^{2}\beta_{1}\gamma_{1}
+b_{1}c_{1}\alpha_{1}^{2}\mu 
\\
& \quad
+2a_{1}c_{1}\alpha_{1}\beta_{1}\mu
+2a_{1}b_{1}\alpha_{1}\gamma_{1}\mu
+a_{1}^{2}\beta_{1}\gamma_{1}\mu)xy -3b(a_{1}b_{1}^{2}c_{1}
+\alpha_{1}\beta_{1}^{2}\gamma_{1} \\
& \quad
+2b_{1}c_{1}\alpha_{1}\beta_{1}\mu
+a_{1}c_{1}\beta_{1}^{2}\mu
+b_{1}^{2}\alpha_{1}\gamma_{1}\mu
+2a_{1}b_{1}\beta_{1}\gamma_{1}\mu)y^{2}+(a^{2}a_{1}^{2} 
-3a_{1}^{2}b c_{1}^{2}\\
& \quad
+2a a_{1}b\alpha_{1}
+b^{2}\alpha_{1}^{2}
-3b\alpha_{1}^{2}\gamma_{1}^{2}
-3b c_{1}^{2}\alpha_{1}^{2}\mu 
-12a_{1}b c_{1}\alpha_{1}\gamma_{1}\mu
-3a_{1}^{2}b\gamma_{1}^{2}\mu)x \\
& \quad+(a^{2}a_{1}b_{1}
-3a_{1}b b_{1}c_{1}^{2}
+ab b_{1}\alpha_{1}
+a a_{1}b\beta_{1}
+b^{2}\alpha_{1}\beta_{1}
-3b\alpha_{1}\beta_{1}\gamma_{1}^{2} \\
& \quad
-3b c_{1}^{2}\alpha_{1}\beta_{1}\mu
-6b b_{1}c_{1}\alpha_{1}\gamma_{1}\mu
-6a_{1}b c_{1}\beta_{1}\gamma_{1}\mu
-3a_{1}b b_{1}\gamma_{1}^{2}\mu)y
\big),
\end{aligned}
\end{equation}
with the first integral
\begin{equation}
\begin{aligned}
H_6^{2}(x,y)=&\frac{a^{2}}{2b}\left(c_{1}+a_{1}x+b_{1}y\right)^{2}
+a\left(c_{1}+a_{1}x+b_{1}y\right)
\left(x\alpha_{1}+y\beta_{1}+\gamma_{1}\right) \nonumber\\
&
+\frac{b}{2}\left(x\alpha_{1}+y\beta_{1}+\gamma_{1}\right)^{2}
+\frac{1}{4}\big(
-\left(c_{1}+a_{1}x+b_{1}y\right)^{4}
-(x\alpha_{1}+y\beta_{1}  \nonumber\\
& +\gamma_{1})^{4}
\big)
-\frac{3}{2}\mu\left(c_{1}+a_{1}x+b_{1}y\right)^{2}
\left(x\alpha_{1}+y\beta_{1}+\gamma_{1}\right)^{2}.
\end{aligned}
\end{equation}

\section{Proof of Theorem $\ref{Teor-Principal}$}
\label{sect:03}
\vspace{0.2cm} 

We consider in $\Sigma^+$ the Hamiltonian nilpotent saddle $(\mathtt{N}_i)$ after an arbitrary affine change of variables,  with its first
integral $H_i(x, y)$; and in $\Sigma^-$,  a  linear differential center  $\mathtt{L}_c$ with its first integral $H_L(x, y)$.

If the discontinuous piecewise differential system $(\mathtt{L}_c)$-$(\mathtt{N}_i)$ for $i=1, \dots, 6$ admits a limit cycle that intersects the nonregular line  $\Sigma$ at two points, $(x,0)$ and $(0, y)$, then these two points must satisfy the following system of equations:
\begin{align}\label{sist_H_Hi}
   e_1=& H_L(x,0) - H_L(0,y)=P_{L}(x,y) = 0, \nonumber \\
e_i=& H_i(x,0) - H_i(0,y)=P_{i}(x,y)  = 0,
\end{align}
where 
\begin{align*}
    P_L (x,y) =& (x+Ay)^2+2(C\,x-By)+\Omega^{2}y^{2}
\end{align*}
and $P_i(x,y)$ are polynomials of degrees two and four respectively, are given in \ref{app1}. 


By Bézout's theorem \cite{coolidge2004treatise}, the number of isolated solutions of this system in $\mathbb{C}^2$ is at most $8$. Consequently, the number of real solutions is also bounded by $8$.

Observe that $(0,0)$ is always a solution, but it does not correspond to a crossing limit cycle. Moreover, each crossing limit cycle is uniquely associated with a solution $(x,y)$ satisfying $x>0$ and $y>0$, since it determines the two intersection points of the orbit with $\Sigma$. Therefore, the number of crossing limit cycles is bounded by the number of admissible solutions $(x,y)$ with $x>0$ and $y>0$, which is at most $7$.




All numerical computations were performed using a default absolute tolerance of $10^{-14}$; decimal expansions were truncated for clarity of presentation.


\subsection*{Limit cycles for the discontinuous piecewise differential system $(\mathtt{L}_c)$-$(\mathtt{N}_1)$}
\begin{figure}[t]
\centering
\includegraphics[scale=0.3]{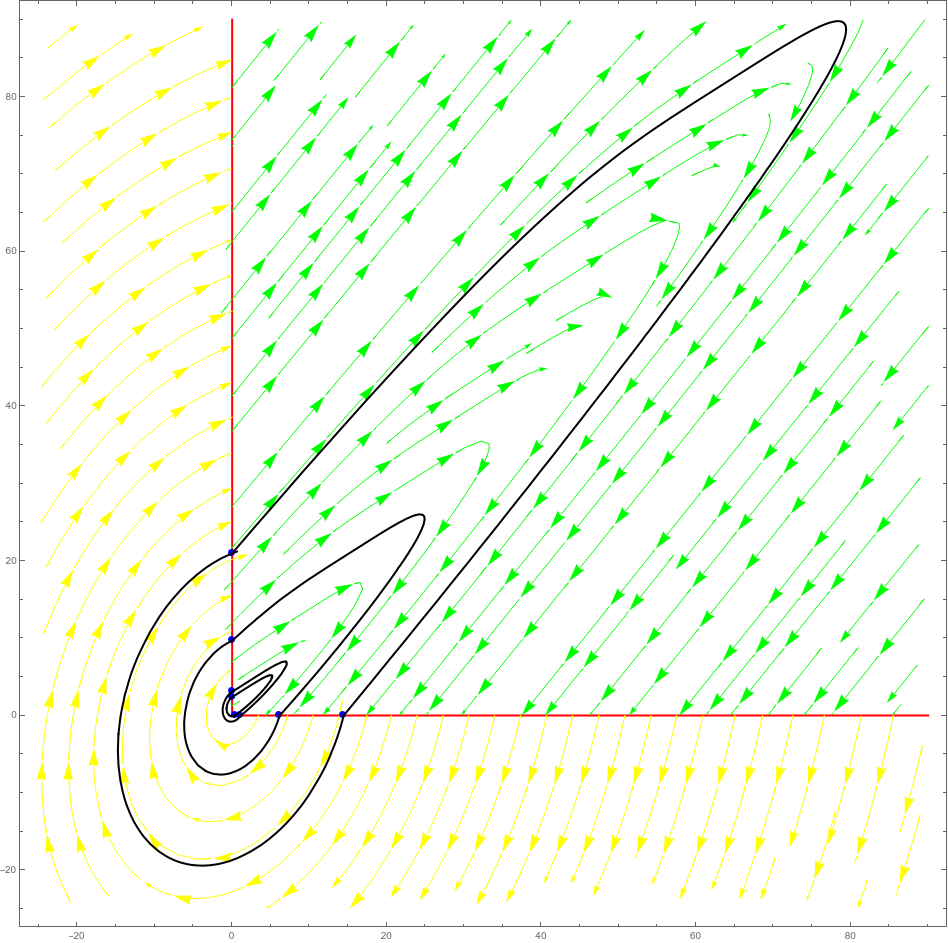}
\caption{The four limit cycles of the discontinuous piecewise differential system \eqref{sist11-Hc-H1}-\eqref{sist12-Hc-H1}.}
\label{fig-Hc-H1}
\end{figure}

Now we give an example of a discontinuous piecewise differential system  $(\mathtt{L}_c)$-$(\mathtt{N}_1)$ separated by the nonregular line  $\Sigma$, having four limit cycles. 
In $\Sigma^-$ we  consider the  linear differential center
\begin{equation}\label{sist11-Hc-H1}
\begin{aligned}
    \dot{x} =& -\frac{6}{5}-\frac{2}{5}\left(x-\frac{y}{5}\right)
               +\frac{49}{50}y, \\[4pt]
    \dot{y} =& -2\left(x-\frac{y}{5}\right),
\end{aligned}
\end{equation}
with the first integral
\begin{equation*}
    H_L(x,y) = \left(x-\frac{y}{5}\right)^2-\frac{6y}{5}+\frac{49y^2}{100}.
\end{equation*}
In $\Sigma^{+}$ we consider the Hamiltonian nilpotent saddle of type $(\mathtt{N}_1)$
\begin{equation}\label{sist12-Hc-H1}
\begin{aligned}
    \dot{x} =& \frac{1}{800}\Big( -3892 - 125x^{3} + 300x^{2}(-2+y) + 4148y
    + 64(-6+y)y^{2}  \\
    &\hspace{2.2cm} - 15x\big(207 + 16(-4+y)y\big) \Big), \\[4pt]
    \dot{y} =& \frac{1}{640}\Big( -125x^{3} + 300x^{2}(-2+y) - 3x\big(683 + 80(-4+y)y\big) \\
    &\hspace{2.2cm} + 4\big(-557 + y(621 + 16(-6+y)y)\big) \Big),
\end{aligned}
\end{equation} 
with the first integral    
\begin{align*}
H_1(x,y) =&
-\frac{2}{5}\left(\frac{3}{10}+\frac{x}{5}-\frac{4y}{5}\right)^2
-\left(\frac{3}{10}+\frac{x}{5}-\frac{4y}{5}\right)
 \left(\frac{4}{5}+\frac{x}{2}-\frac{2y}{5}\right) \notag\\
&\quad
-\frac{5}{8}\left(\frac{4}{5}+\frac{x}{2}-\frac{2y}{5}\right)^2
-\frac{1}{4}\left(\frac{4}{5}+\frac{x}{2}-\frac{2y}{5}\right)^4 .
\end{align*}
Solving the system~\eqref{sist_H_Hi}, we obtain four real pairs of solutions
$(p_k,q_k)$, $k=1,\dots,4$, where $p_k=(x_k,0)$ and $q_k=(0,y_k)$. More precisely,
\[
\begin{aligned}
p_1 &\approx (0.387552,0), & q_1 &\approx (0,2.38307),
& p_2 &\approx (1.13899,0), & q_2 &\approx (0,3.06322),\\
p_3 &\approx (6.15242,0), & q_3 &\approx (0,9.65856),
& p_4 &\approx (14.4234,0), & q_4 &\approx (0,20.9765).
\end{aligned}
\]
that provide the crossing limit cycles of discontinuous piecewise differential system \eqref{sist11-Hc-H1}–\eqref{sist12-Hc-H1} shown in Figure \ref{fig-Hc-H1}. 


\subsection*{Limit cycles for the discontinuous piecewise differential system $(\mathtt{L}_c)$-$(\mathtt{N}_2)$}

\begin{figure}[t]
\begin{center}
\includegraphics[scale=0.3]{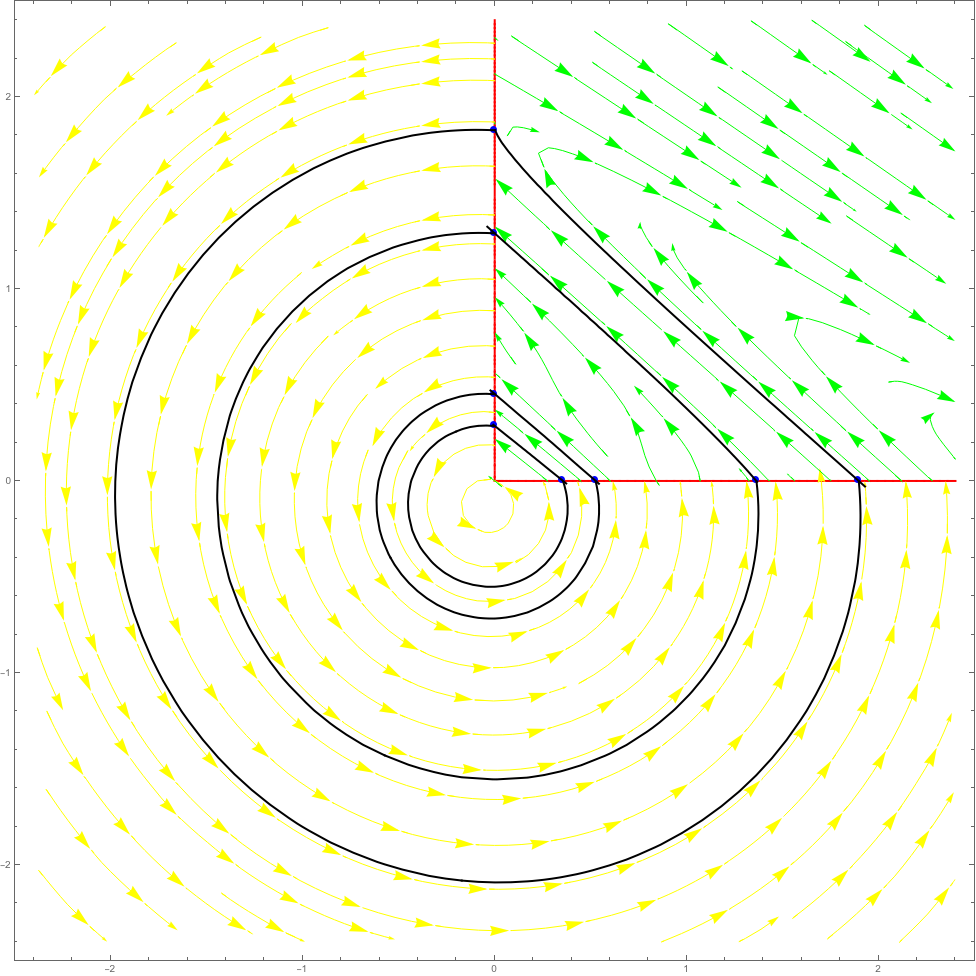}
\caption{The four limit cycles of the discontinuous piecewise differential system \eqref{sist21-Hc-H2}-\eqref{sist22-Hc-H2}.  }
\label{fig-Hc-H2}
\end{center}
\end{figure}

Now we give an example of a discontinuous piecewise differential system  $(\mathtt{L}_c)$-$(\mathtt{N}_2)$ separated by the nonregular line  $\Sigma$, having four limit cycles. 
In $\Sigma^-$ we  consider the  linear differential center
\begin{equation}\label{sist21-Hc-H2}
\begin{aligned}
    \dot{x} =& -\frac{13}{50}
               -\frac{3}{50}\left(x+\frac{3}{100}y\right)
               -\frac{9801}{5000}y, \\[4pt]
    \dot{y} =& \frac{2}{25}
               +2\left(x+\frac{3}{100}y\right),
\end{aligned}
\end{equation}
with the first integral
\begin{equation*}
    H_L(x,y) =-\left(x+\frac{3}{100}y\right)^2 - 2\left(\frac{x}{25}+\frac{13}{100}y\right) - \frac{9801}{10000}y^{2}.
\end{equation*}
In $\Sigma^{+}$ we consider the Hamiltonian nilpotent saddle of type $(\mathtt{N}_2)$
\begin{equation}\label{sist22-Hc-H2}
\begin{aligned}
\dot{x} =& \frac{1}{119910000}\Big( -1072596688 + 2477433732\,x +84672x^{2}\big(-21011+4380x\big) \\
&\quad
 + 4\big(818403452+4473x(-267317+85128x)\big)y +317583(-10020 + \\
&\quad
6523x)y^{2} + 932000244\,y^{3} \Big), \\[6pt]
\dot{y} =& \frac{1}{39970000}\Big( 275724936 -48x\big(12792695+2016x(-4521+928x)\big) - 4 \\
&\quad
\big(206452811+14112x(-21011+6570x)\big)y + 2982(267317- 170256 \\
&\quad
x)y^{2} - 230177101\,y^{3} \Big),
\end{aligned}
\end{equation}
with the first integral    
\begin{align*}
 H_2(x,y) =&\frac{7}{50}\left(\frac{27}{50}-\frac{87}{100}x-\frac{93}{100}y\right)^2
+\frac{2}{25}\left(\frac{27}{50}-\frac{87}{100}x-\frac{93}{100}y\right) \\
&\quad
 \left(\frac{47}{50}-\frac{12}{25}x-\frac{71}{100}y\right) +\frac{2}{175}\left(\frac{47}{50}-\frac{12}{25}x-\frac{71}{100}y\right)^2 - \\
&\quad
\left(\frac{27}{50}- \frac{87}{100}x-\frac{93}{100}y\right)
 \left(\frac{47}{50}-\frac{12}{25}x-\frac{71}{100}y\right)^3.
\end{align*}
Solving the system~\eqref{sist_H_Hi}, we obtain four real pairs of solutions
$(p_k,q_k)$, $k=1,\dots,4$, where $p_k=(x_k,0)$ and $q_k=(0,y_k)$. More precisely,
\[
\begin{aligned}
p_1 &\approx (0.355545,0), & q_1 &\approx (0,0.286309),
& p_2 &\approx (0.525244,0), & q_2 &\approx (0,0.451964),\\
p_3 &\approx (1.36335,0), & q_3 &\approx (0,1.28996),
& p_4 &\approx (1.89636,0), & q_4 &\approx (0,1.82657).
\end{aligned}
\]
that provide the crossing limit cycles of discontinuous piecewise differential system \eqref{sist21-Hc-H2}–\eqref{sist22-Hc-H2} shown in Figure \ref{fig-Hc-H2}. 

 
\subsection*{Limit cycles for the discontinuous piecewise differential system $(\mathtt{L}_c)$-$(\mathtt{N}_3^{1})$}
\begin{figure}[t]
\begin{center}
\includegraphics[width=0.45\textwidth, height=0.3 \textheight]{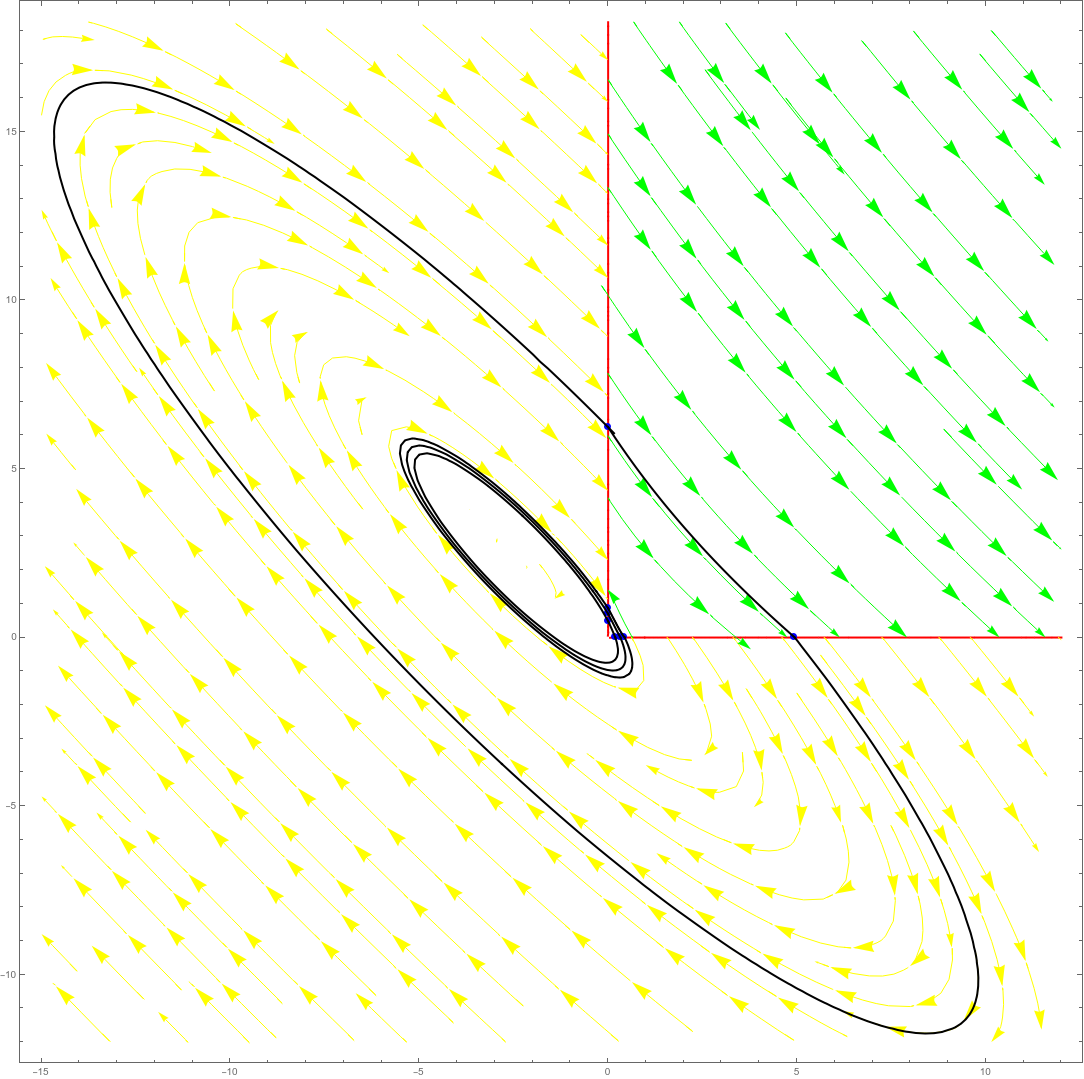}
\includegraphics[width=0.45\textwidth, height=0.3 \textheight]{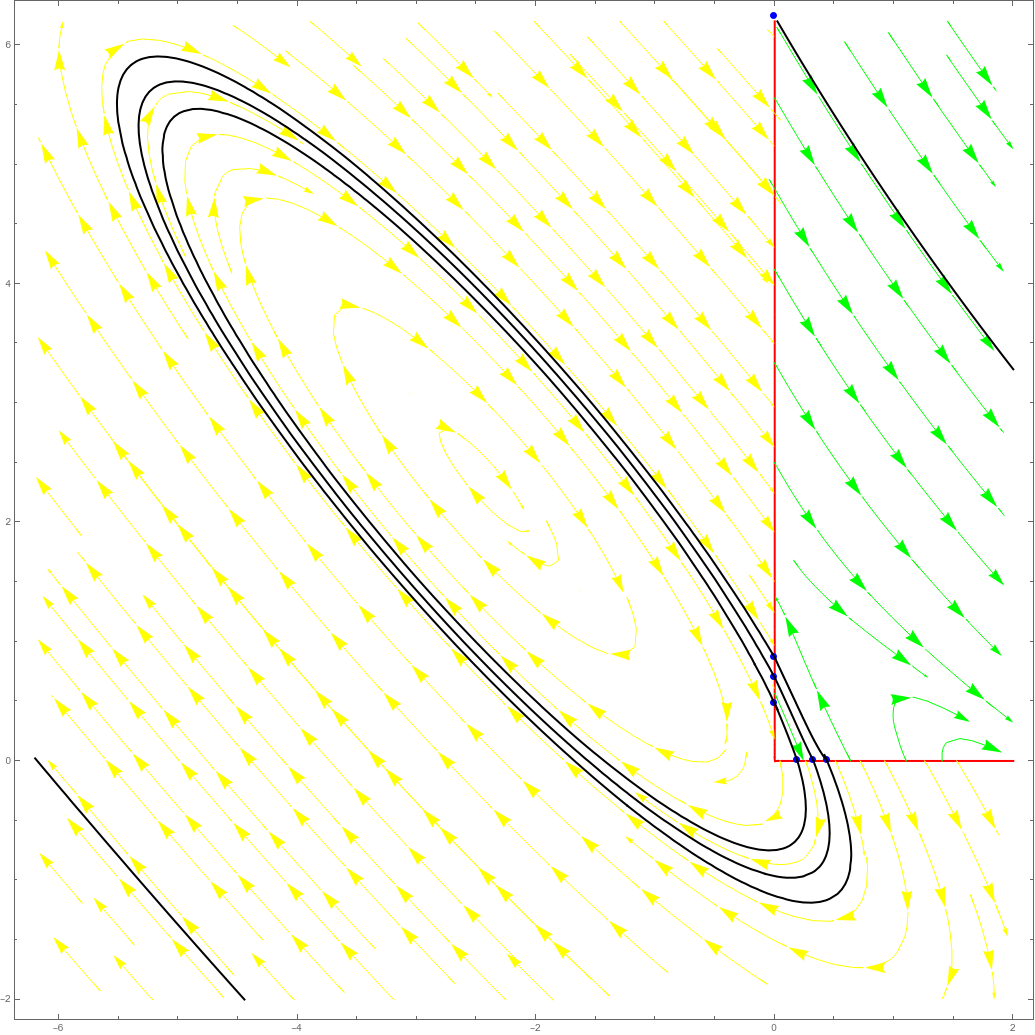}
\caption{The four limit cycles of the discontinuous piecewise differential system \eqref{eq:Xh-11}–\eqref{eq:Yf-12}. Right: zoom around the origin.}
\label{fig-Hc-H31}
\end{center}
\end{figure}

Now we give an example of a discontinuous piecewise differential system  $(\mathtt{L}_c)$-$(\mathtt{N}_3^{1})$ separated by the nonregular line  $\Sigma$, having four limit cycles. 
In $\Sigma^-$ we  consider the  linear differential center
\begin{equation}\label{eq:Xh-11}
\begin{aligned}
\dot{x} =& \frac{1}{5}
+\frac{77}{50}\left(x+\frac{77}{100}y\right)
+\frac{8}{25}y, \\[4pt]
\dot{y} =& -\frac{31}{25}
-2\left(x+\frac{77}{100}y\right),
\end{aligned}
\end{equation}
with the first integral
\begin{align*}
H_{L}(x,y) =& 2\left(\frac{31}{50}x+\frac{1}{10}y\right) + \left(x+\frac{77}{100}y\right)^2 + \frac{4}{25}y^2.
\end{align*}
In $\Sigma^{+}$ we consider the Hamiltonian nilpotent saddle of type $(\mathtt{N}_3^{1})$
\begin{equation}\label{eq:Yf-12}
\begin{aligned}
\dot{x} =& \frac{1}{4130000} \Big( -405967+5474739\,y+9\Big(778786\,x +13x^{2}\big(-146848 \\
& +73877x\big)+4x\big(-728887+615690x\big)y+57\big(-16337+31442x\big)y^{2} \\
& +357029\,y^{3}\big)\Big),\\[6pt]
\dot{y} =& \frac{1}{12390000}\Big( -27263925\,x^{3}+1053\big(55627-73877y\big)x^{2} -24x\big(716899 \\
& +9y(-477256+307845y)\big) -2\big(862465+27y\big(389393+y(-728887 \\
&\qquad
+298699y)\big)\big) \Big),
\end{aligned}
\end{equation}
with the first integral    
\begin{align*}
H_3^{1}(x,y) &=\frac{73}{200}\left(\frac{43}{100}-\frac{24}{25}x-\frac{9}{20}y\right)^2
-\frac{3}{2}\left(\frac{43}{100}-\frac{24}{25}x-\frac{9}{20}y\right)^2 \nonumber\\
 &
 \left(-\frac{39}{100}+\frac{39}{100}x+\frac{57}{100}y\right)^2 +\frac{1}{4}\left(-\frac{39}{100}+\frac{39}{100}x+\frac{57}{100}y\right)^4.
\end{align*}
Solving the system~\eqref{sist_H_Hi}, we obtain four real pairs of solutions
$(p_k,q_k)$, $k=1,\dots,4$, where $p_k=(x_k,0)$ and $q_k=(0,y_k)$. More precisely,
\[
\begin{aligned}
p_1 &\approx (0.190098,0), & q_1 &\approx (0,0.482586),
& p_2 &\approx (0.325214,0), & q_2 &\approx (0,0.700087),\\
p_3 &\approx (0.439849,0), & q_3 &\approx (0,0.86669),
& p_4 &\approx (4.94215,0), & q_4 &\approx (0,6.23885).
\end{aligned}
\]
that provide the crossing limit cycles of discontinuous piecewise differential system \eqref{eq:Xh-11}–\eqref{eq:Yf-12} shown in Figure \ref{fig-Hc-H31}. 
 
\subsection*{Limit cycles for the discontinuous piecewise differential system $(\mathtt{L}_c)$-$(\mathtt{N}_3^{2})$}
Now we give an example of a discontinuous piecewise differential system  $(\mathtt{L}_c)$-$(\mathtt{N}_3^{2})$ separated by the nonregular line  $\Sigma$, having four limit cycles. 
In $\Sigma^-$ we  consider the  linear differential center
\begin{equation}\label{eq:Xh-21}
\begin{aligned}
\dot{x} =& \frac{3}{5}+x+\frac{37}{25}y, \\[4pt]
\dot{y} =& \frac{4}{5}-2\left(x+\frac{y}{2}\right),
\end{aligned}
\end{equation}
with the first integral
\begin{equation*}
    H_L(x,y) =2\left(-\frac{2}{5}x+\frac{3}{10}y\right) + \left(x+\frac{y}{2}\right)^2 + \frac{49}{100}y^2.
\end{equation*}
In $\Sigma^{+}$ we consider the Hamiltonian nilpotent saddle of type $(\mathtt{N}_3^{2})$
\begin{equation}\label{eq:Yf-22}
\begin{aligned}
\dot{x} =& \frac{1}{1500}\Big(
-22688+55x^{3}+6x^{2}(307-50y)-24x\big(-615+2y(13+5y)\big) \\
&\quad
+y\big(31560+y(-6618+485y)\big)\Big), \\[6pt]
\dot{y} =& \frac{1}{1500}\Big( -215x^{3}-3x^{2}(676+55y)+12x\big(-510+y(-307+25y)\big) \\
&\quad
+8\big(1564+y(-1845+y(39+10y))\big)\Big),
\end{aligned}
\end{equation}
with the first integral    
\begin{equation*}
\begin{aligned}
 H_3^{2}(x,y)=&
-\frac{8}{15}\left(1-\frac{3}{5}x-\frac{9}{10}y\right)^2
-\frac{4}{5}\left(1-\frac{3}{5}x-\frac{9}{10}y\right)
 \left(-\frac{4}{5}-\frac{x}{10}+\frac{y}{10}\right) \\
& -\frac{3}{10}\left(-\frac{4}{5}-\frac{x}{10}+\frac{y}{10}\right)^2
-\frac{3}{2}\left(1-\frac{3}{5}x-\frac{9}{10}y\right)^2
 \left(-\frac{4}{5}-\frac{x}{10}+\frac{y}{10}\right)^2 \\
& +\frac{1}{4}\left(-\frac{4}{5}-\frac{x}{10}+\frac{y}{10}\right)^4.
\end{aligned}
\end{equation*}

\begin{figure}[t]
\begin{center}
\includegraphics[width=0.45\textwidth, height=0.3 \textheight]{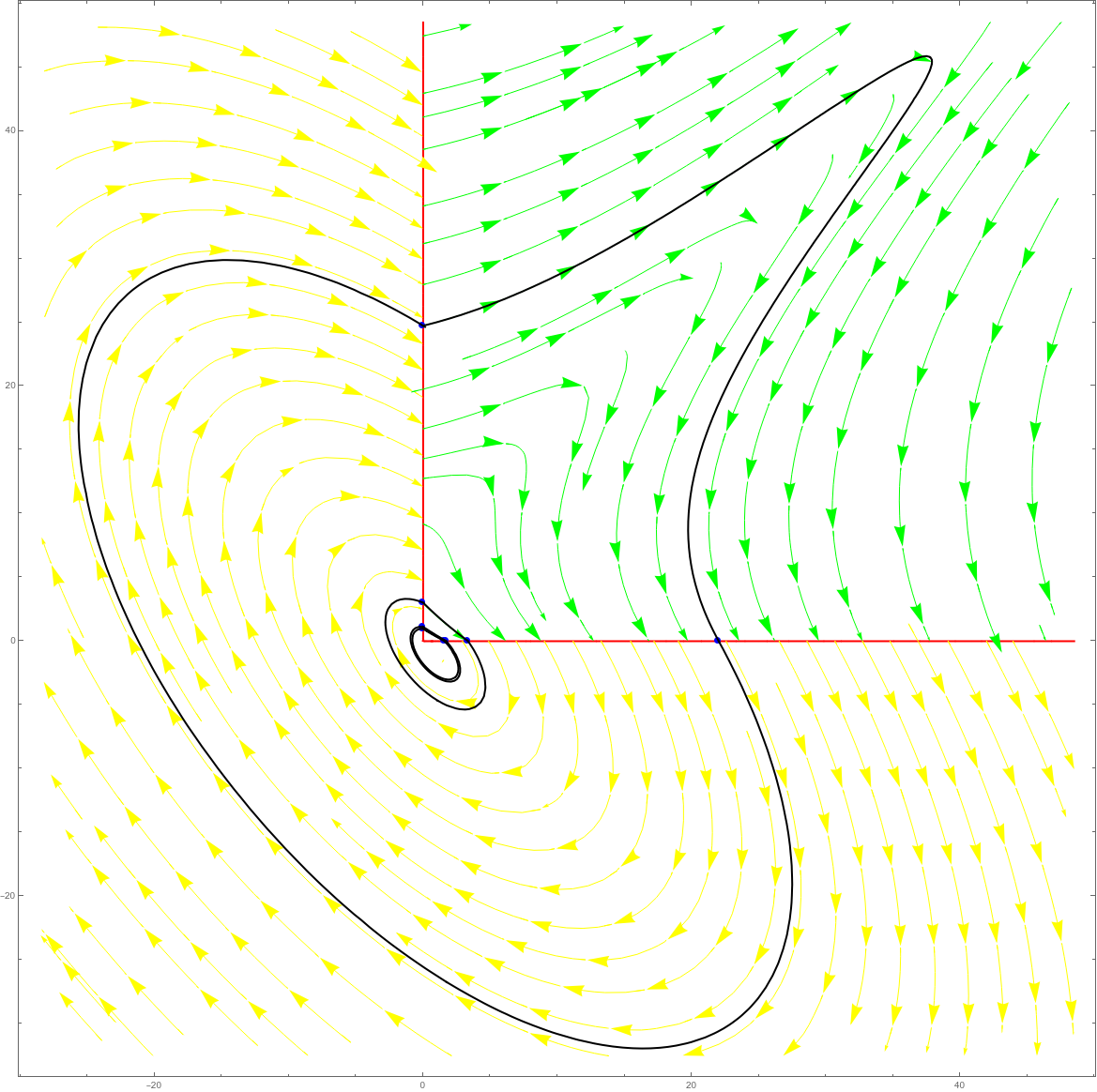}
\includegraphics[width=0.45\textwidth, height=0.3 \textheight]{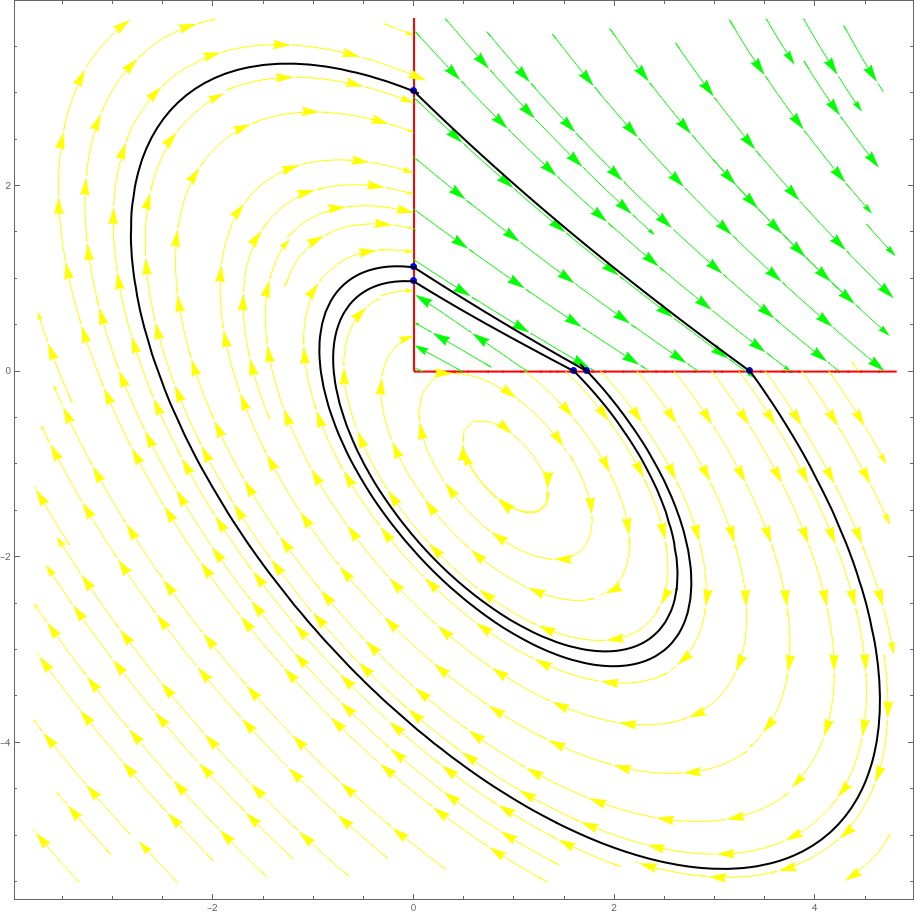}
\caption{The four limit cycles of the discontinuous piecewise differential system \eqref{eq:Xh-21}–\eqref{eq:Yf-22}. Right: zoom around the origin.}
\label{fig-Hc-H32}
\end{center}
\end{figure}

Solving the system~\eqref{sist_H_Hi}, we obtain four real pairs of solutions
$(p_k,q_k)$, $k=1,\dots,4$, where $p_k=(x_k,0)$ and $q_k=(0,y_k)$. More precisely,
\[
\begin{aligned}
p_1 &\approx (1.60038,0), & q_1 &\approx (0,0.971298),
& p_2 &\approx (1.72908,0), & q_2 &\approx (0,1.12275),\\
p_3 &\approx (3.35256,0), & q_3 &\approx (0,3.01931),
& p_4 &\approx (22.0218,0), & q_4 &\approx (0,24.7284).
\end{aligned}
\]

that provide the crossing limit cycles of discontinuous piecewise differential system \eqref{eq:Xh-21}–\eqref{eq:Yf-22} shown in Figure \ref{fig-Hc-H32}. 
 
 \subsection*{Limit cycles for the discontinuous piecewise differential system $(\mathtt{L}_c)$-$(\mathtt{N}_4^{1})$ }
\begin{figure}[t]
\begin{center}
\includegraphics[scale=0.3]{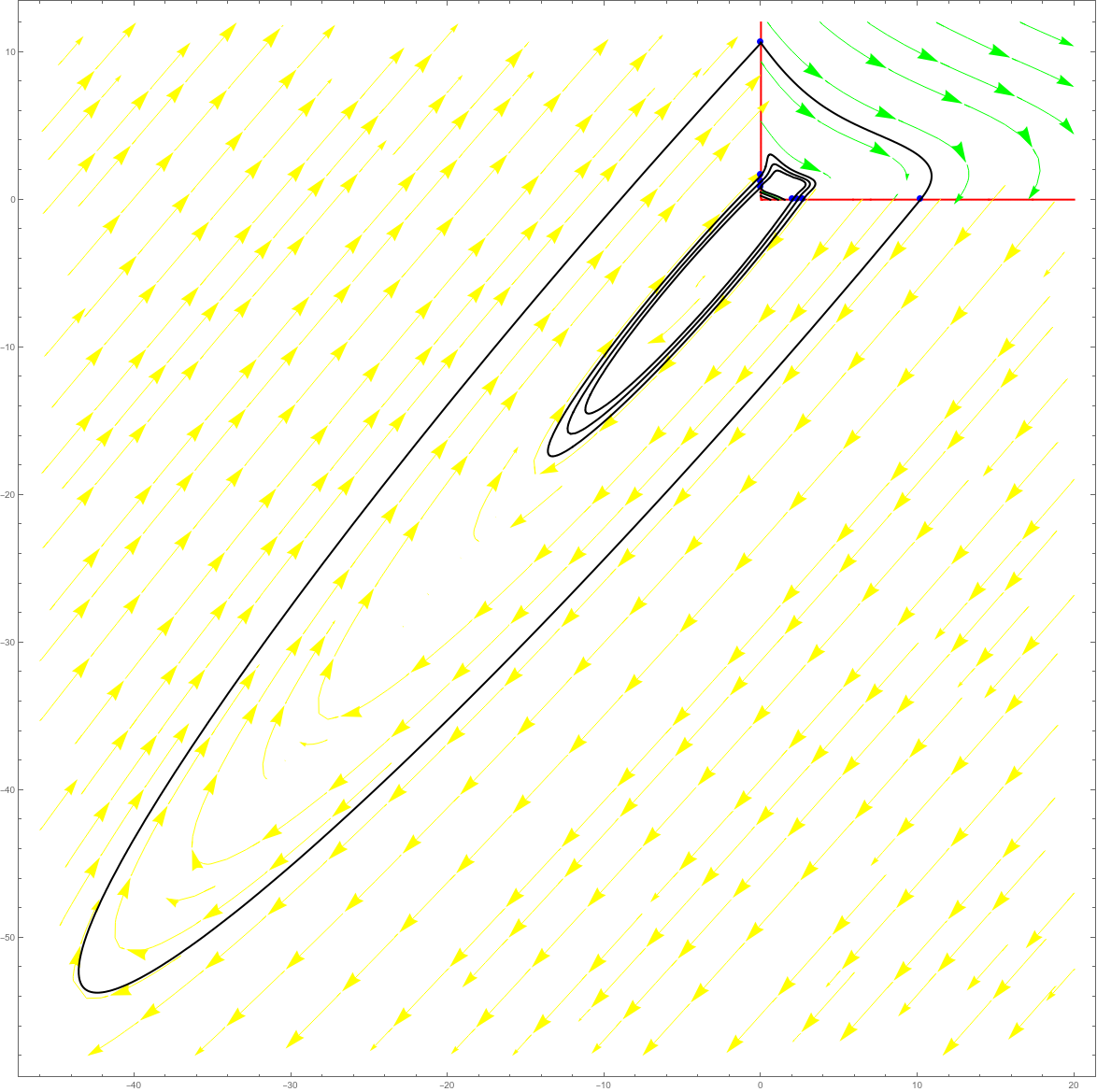}
\caption{The four limit cycles of the discontinuous piecewise differential system \eqref{eq:Xh-41}–\eqref{eq:Yf-41}.}
\label{fig-Hc-H41}
\end{center}
\end{figure}

Now we give an example of a discontinuous piecewise differential system  $(\mathtt{L}_c)$-$(\mathtt{N}_4^{1})$ separated by the nonregular line  $\Sigma$, having four limit cycles. 
In $\Sigma^-$ we  consider the  linear differential center

\begin{equation}\label{eq:Xh-41}
\begin{aligned}
\dot{x} =& \frac{7}{5}-\frac{8}{5}\left(x-\frac{4}{5}y\right)+\frac{2}{25}y, \\[4pt]
\dot{y} =& \frac{6}{5}-2\left(x-\frac{4}{5}y\right),
\end{aligned}
\end{equation}
with the first integral
\begin{equation*}
    H_L(x,y) =\left(x-\frac{4}{5}y\right)^2+2\left(-\frac{3}{5}x+\frac{7}{10}y\right)+\frac{1}{25}y^2.
\end{equation*}
In $\Sigma^{+}$ we consider the Hamiltonian nilpotent saddle of type $(\mathtt{N}_4^{1})$
\begin{equation}\label{eq:Yf-41}
\begin{aligned}
\dot{x} =& \frac{1}{2100}\Big(-8013+x(6859-2x(621+80x))+14302y \\
&\quad
+12x(-927+173x)y+12(-513+233x)y^{2}+776y^{3}\Big), \\[6pt]
\dot{y} =& \frac{1}{4200}\Big(6927-775x^{3}+30x^{2}(81+32y) \\
&\quad
+x\big(-5798+24(207-173y)y\big)-2y(6859-5562y+932y^{2})\Big),
\end{aligned}
\end{equation}
with the first integral    
\begin{equation*}
\begin{aligned}
 H_{41}(x,y) =&-\frac{1}{4}\left(-\frac{9}{10}+\frac{x}{2}+\frac{y}{5}\right)^4
-\frac{1}{10}\left(-\frac{3}{5}-\frac{x}{10}+\frac{4}{5}y\right)^2 \\
&-\frac{3}{2}\left(-\frac{9}{10}+\frac{x}{2}+\frac{y}{5}\right)^2
 \left(-\frac{3}{5}-\frac{x}{10}+\frac{4}{5}y\right)^2.
\end{aligned}
\end{equation*} 
Solving the system~\eqref{sist_H_Hi}, we obtain four real pairs of solutions
$(p_k,q_k)$, $k=1,\dots,4$, where $p_k=(x_k,0)$ and $q_k=(0,y_k)$. More precisely,
\[
\begin{aligned}
p_1 &\approx (2.02448,0), & q_1 &\approx (0,0.845234),
& p_2 &\approx (2.35986,0), & q_2 &\approx (0,1.22555),\\
p_3 &\approx (2.70908,0), & q_3 &\approx (0,1.62987),
& p_4 &\approx (10.1815,0), & q_4 &\approx (0,10.6126).
\end{aligned}
\]
that provide the crossing limit cycles of discontinuous piecewise differential system \eqref{eq:Xh-41}–\eqref{eq:Yf-41} shown in Figure \ref{fig-Hc-H41}. 
 
 \subsection*{Limit cycles for the discontinuous piecewise differential system $(\mathtt{L}_c)$-$(\mathtt{N}_4^{2})$ }
\begin{figure}[t]
\begin{center}
\includegraphics[scale=0.3]{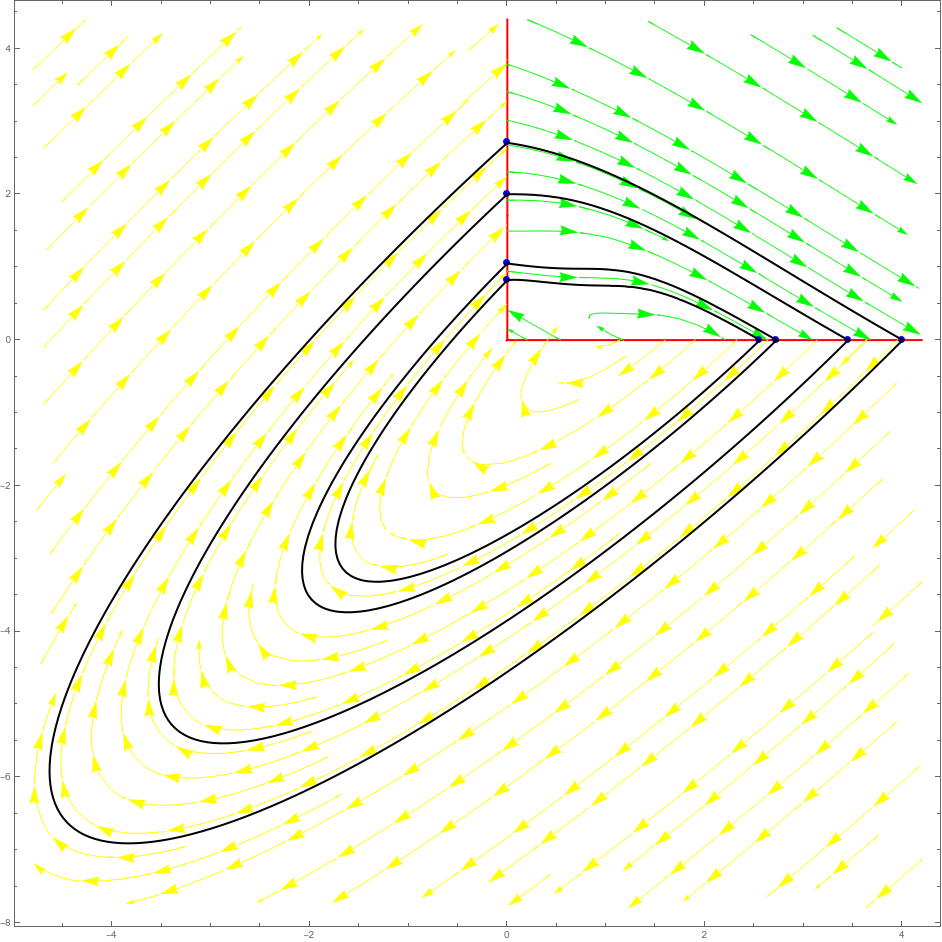}
\caption{The four limit cycles of the discontinuous piecewise differential system \eqref{eq:Xh-42}–\eqref{eq:Yf-42}.}
\label{fig-Hc-H42}
\end{center}
\end{figure}

Now we give an example of a discontinuous piecewise differential system  $(\mathtt{L}_c)$-$(\mathtt{N}_4^{2})$  separated by the nonregular line  $\Sigma$, having four limit cycles. 
In $\Sigma^-$ we  consider the  linear differential center
\begin{equation}\label{eq:Xh-42}
\begin{aligned}
\dot{x} =& \frac{6}{5} - \frac{7}{5}\left(x-\frac{7}{10}y\right) + \frac{8}{25}y, \\[4pt]
\dot{y} =& 2 - 2\left(x-\frac{7}{10}y\right),
\end{aligned}
\end{equation}
 with the first integral
\begin{equation*}
H_L(x,y) =\left(x-\frac{7}{10}y\right)^2 + 2\left(-x+\frac{3}{5}y\right) + \frac{4}{25}y^2.
\end{equation*}
In $\Sigma^{+}$ we consider the Hamiltonian nilpotent saddle of type $(\mathtt{N}_4^{2})$

\begin{equation}\label{eq:Yf-42}
\begin{aligned}
\dot{x} =& \frac{1}{3200}\Big( -16108+x(28511+8x(-2127+424x))+34726y \\
&\quad
+96x(-413+124x)y+96(-229+136x)y^{2}+5632y^{3}\Big), \\[6pt]
\dot{y} =& \frac{1}{6400}\Big( 25276+x(-45083+8(3351-664x)x)-57022y \\
&\quad
+96(709-212x)xy+96(413-248x)y^{2}-8704y^{3}\Big),
\end{aligned}
\end{equation}
with the first integral    
\begin{equation*}
\begin{aligned}
 H_{4}^{2}(x,y) =&-\frac{1}{5}\left(\frac{7}{10}-\frac{2}{5}x-\frac{4}{5}y\right)^2
-\frac{1}{4}\left(\frac{7}{10}-\frac{2}{5}x-\frac{4}{5}y\right)^4 \\
&-\frac{3}{10}\left(\frac{7}{10}-\frac{2}{5}x-\frac{4}{5}y\right)
 \left(-\frac{4}{5}+\frac{x}{2}+\frac{y}{5}\right)
-\frac{9}{80}\left(-\frac{4}{5}+\frac{x}{2}+\frac{y}{5}\right)^2 \\
&-\frac{3}{2}\left(\frac{7}{10}-\frac{2}{5}x-\frac{4}{5}y\right)^2
 \left(-\frac{4}{5}+\frac{x}{2}+\frac{y}{5}\right)^2 .
\end{aligned}
\end{equation*} 
Solving the system~\eqref{sist_H_Hi}, we obtain four real pairs of solutions
$(p_k,q_k)$, $k=1,\dots,4$, where $p_k=(x_k,0)$ and $q_k=(0,y_k)$. More precisely,
\[
\begin{aligned}
p_1 &\approx (2.55713,0), & q_1 &\approx (0,0.821581),
& p_2 &\approx (2.72657,0), & q_2 &\approx (0,1.05173),\\
p_3 &\approx (3.4514,0), & q_3 &\approx (0,2.00246),
& p_4 &\approx (4.00261,0), & q_4 &\approx (0,2.70789).
\end{aligned}
\]
that provide the crossing limit cycles of discontinuous piecewise differential system \eqref{eq:Xh-42}–\eqref{eq:Yf-42} shown in Figure \ref{fig-Hc-H42}. 
 
 \subsection*{Limit cycles for the discontinuous piecewise differential system $(\mathtt{L}_c)$-$(\mathtt{N}_5^{1})$}
\begin{figure}[t]
\begin{center}
\includegraphics[scale=0.3]{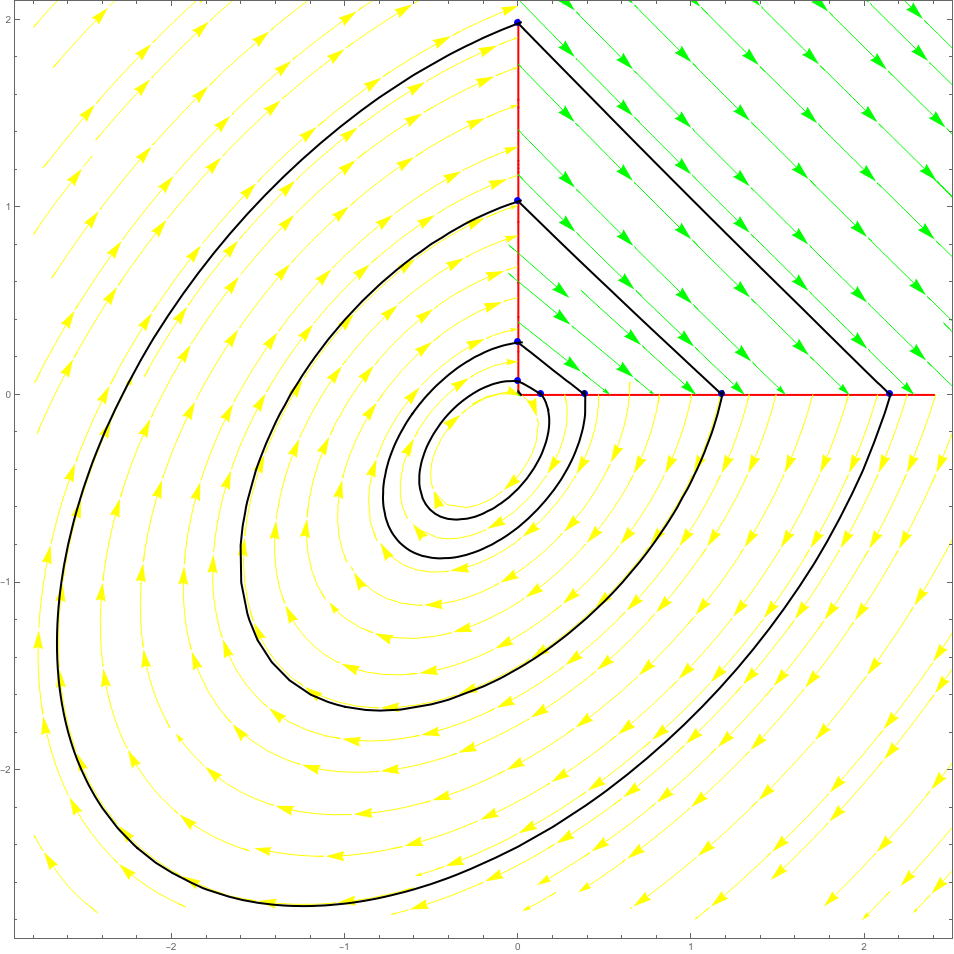}
\caption{The four limit cycles of the discontinuous piecewise differential system \eqref{eq:Xh-51}–\eqref{eq:Yf-51}.}
\label{fig-Hc-H51}
\end{center}
\end{figure}

Now we give an example of a discontinuous piecewise differential system  $(\mathtt{L}_c)$-$(\mathtt{N}_5^{1})$ separated by the nonregular line  $\Sigma$, having four limit cycles. 
In $\Sigma^-$ we  consider the  linear differential center

\begin{equation}\label{eq:Xh-51}
\begin{aligned}
\dot{x} =& \frac{11}{25} - \frac{43}{50}\left(x-\frac{43}{100}y\right) + \frac{1058}{625}y, \\[4pt]
\dot{y} =& -\frac{7}{50}-2\left(x-\frac{43}{100}y\right),
\end{aligned}
\end{equation}
with the first integral
\begin{equation*}
    H_L(x,y) =\left(x-\frac{43}{100}y\right)^2 + 2\left(\frac{7}{100}x+\frac{11}{50}y\right)
+\frac{529}{625}y^2.
\end{equation*}
In $\Sigma^{+}$ we consider the Hamiltonian nilpotent saddle of type $(\mathtt{N}_5^{1})$

\begin{equation}\label{eq:Yf-51}
\begin{aligned}
\dot{x} =& \frac{1}{140000000}\Big( -17828080+398806366x^{3}+6x^{2}(-97025863+222074762y) \\
& +2y\big(215022709+68y(-5722053+3962714y)\big) + x\big(354597367+24y \\
& (-55989761+61426967y)\big)\Big), \\[6pt]
\dot{y} =& \frac{1}{280000000}\Big( 38035040+x\big(-590906221-918x(-1102793+776391x)\big) \\
& -709194734y+12(194051726-199403183x)xy + 24(55989761- \\
& 111037381x)y^{2}-982831472y^{3}\Big),
\end{aligned}
\end{equation}
with the first integral    
\begin{equation*}
    \begin{aligned}
 H_{51}(x,y) =&
\frac{1}{4}\big(
-\left(\frac{1}{20}+\frac{19}{50}x+\frac{13}{25}y\right)^4
+\left(-\frac{7}{25}+\frac{59}{100}x+\frac{33}{50}y\right)^4
\big) \\[4pt]
& +\frac{13}{100}\left(\frac{1}{20}+\frac{19}{50}x+\frac{13}{25}y\right)^2 + \frac{21}{100}\left(\frac{1}{20}+\frac{19}{50}x+\frac{13}{25}y\right)^2 \\[4pt]
& \left(-\frac{7}{25}+\frac{59}{100}x+\frac{33}{50}y\right)^2 .
    \end{aligned}
\end{equation*} 
Solving the system~\eqref{sist_H_Hi}, we obtain four real pairs of solutions
$(p_k,q_k)$, $k=1,\dots,4$, where $p_k=(x_k,0)$ and $q_k=(0,y_k)$. More precisely,
\[
\begin{aligned}
p_1 &\approx (0.135002,0), & q_1 &\approx (0,0.072169),
& p_2 &\approx (0.385675,0), & q_2 &\approx (0,0.278707),\\
p_3 &\approx (1.17787,0), & q_3 &\approx (0,1.03194),
& p_4 &\approx (2.14886,0), & q_4 &\approx (0,1.98091).
\end{aligned}
\]
that provide the crossing limit cycles of discontinuous piecewise differential system \eqref{eq:Xh-51}–\eqref{eq:Yf-51} shown in Figure \ref{fig-Hc-H51}. 
 
 
 \subsection*{Limit cycles for the discontinuous piecewise differential system $(\mathtt{L}_c)$-$(\mathtt{N}_5^{2})$}
\begin{figure}[t]
\begin{center}
\includegraphics[scale=0.3]{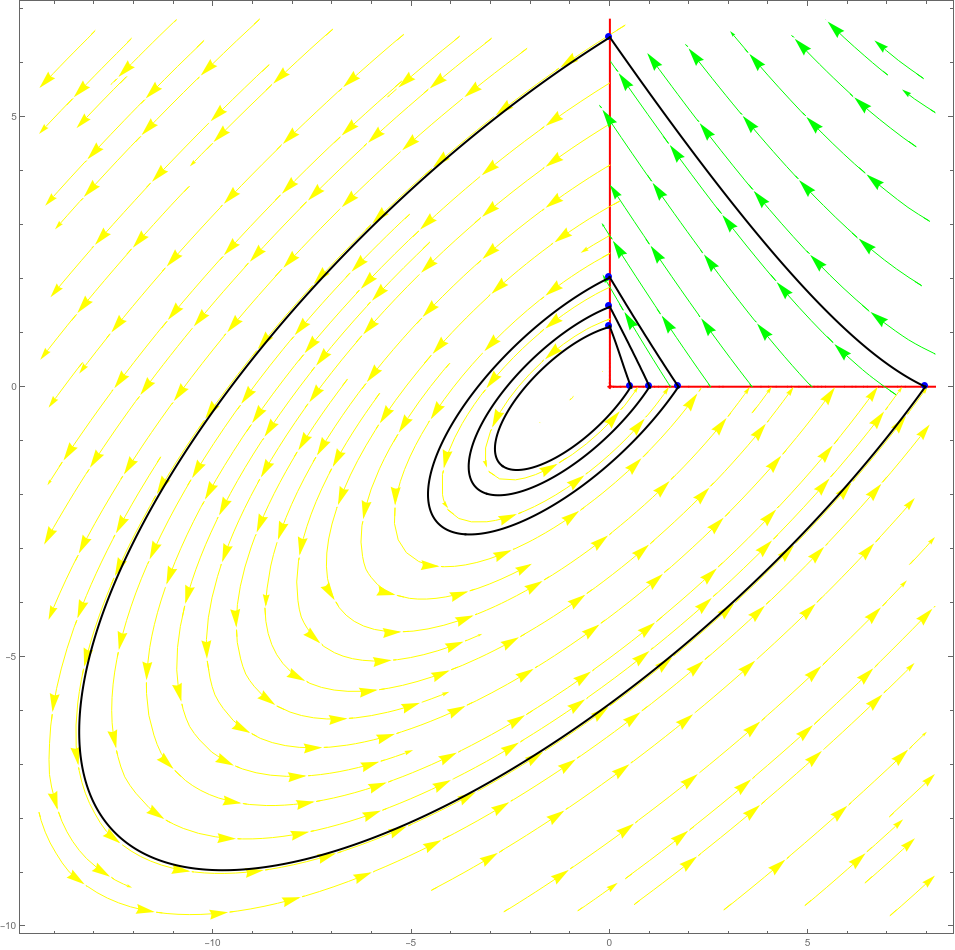}
\caption{The four limit cycles of the discontinuous piecewise differential system \eqref{eq:Xh-52}–\eqref{eq:Yf-52}.}
\label{fig-Hc-H52}
\end{center}
\end{figure}

Now we give an example of a discontinuous piecewise differential system  $(\mathtt{L}_c)$-$(\mathtt{N}_5^{2})$ separated by the nonregular line  $\Sigma$, having four limit cycles. 
In $\Sigma^-$ we  consider the  linear differential center
\begin{equation}\label{eq:Xh-52}
\begin{aligned}
\dot{x} =& \frac{6}{5}+2(x-y)-2y, \\[4pt]
\dot{y} =& \frac{8}{5}+2(x-y),
\end{aligned}
\end{equation}
with the first integral
\begin{equation*}
    H_L(x,y) =-(x-y)^2-2\left(\frac{4}{5}x-\frac{3}{5}y\right)-y^2.
\end{equation*}
In $\Sigma^{+}$ we consider the Hamiltonian nilpotent saddle of type $(\mathtt{N}_5^{2})$
\begin{equation}\label{eq:Yf-52}
\begin{aligned}
\dot{x} =& \frac{1}{1500}\Big( -9995+5x\big(4661+9x(-352+17x)\big) +3\big(5816+ \\
&\quad 9(-864+x)x\big)y  -81(106+17x)y^{2}-729y^{3} \Big), \\[6pt]
\dot{y} =& \frac{1}{9000}\Big( 90395+5x\big(-38291+6(4323-388x)x\big) -30\big(4661 \\
&\quad+9x(-704+51x)\big) y-162(-432+x)y^{2}+2754y^{3}\Big),
\end{aligned}
\end{equation}
with the first integral   
\begin{equation*}
\begin{aligned}
 H_{5}^{2}(x,y) =&
\frac{1}{4}\big( -\left(-\frac{1}{10}-\frac{7}{10}x-\frac{3}{5}y\right)^4
+\left(1-\frac{x}{2}-\frac{3}{10}y\right)^4 \big) +\frac{5}{16} \\[4pt]
& \left(-\frac{1}{10}-\frac{7}{10}x-\frac{3}{5}y\right)^2
+\frac{1}{2}\left(-\frac{1}{10}-\frac{7}{10}x-\frac{3}{5}y\right)
 \left(1-\frac{x}{2}-\frac{3}{10}y\right) \\[4pt]
& +\frac{1}{5}\left(1-\frac{x}{2}-\frac{3}{10}y\right)^2
+\frac{3}{5}\left(-\frac{1}{10}-\frac{7}{10}x-\frac{3}{5}y\right)^2
 \left(1-\frac{x}{2}-\frac{3}{10}y\right)^2 .
\end{aligned}
\end{equation*} 
Solving the system~\eqref{sist_H_Hi}, we obtain four real pairs of solutions
$(p_k,q_k)$, $k=1,\dots,4$, where $p_k=(x_k,0)$ and $q_k=(0,y_k)$. More precisely,
\[
\begin{aligned}
p_1 &\approx (0.52839,0), & q_1 &\approx (0,1.10766),
& p_2 &\approx (1.00057,0), & q_2 &\approx (0,1.47942),\\
p_3 &\approx (1.72915,0), & q_3 &\approx (0,2.02288),
& p_4 &\approx (7.95553,0), & q_4 &\approx (0,6.47249).
\end{aligned}
\]
that provide the crossing limit cycles of discontinuous piecewise differential system \eqref{eq:Xh-52}–\eqref{eq:Yf-52} shown in Figure \ref{fig-Hc-H52}. 

 
\subsection*{Limit cycles for the discontinuous piecewise differential system $(\mathtt{L}_c)$-$(\mathtt{N}_6^{1})$}
\begin{figure}[t]
\begin{center}
\includegraphics[scale=0.3]{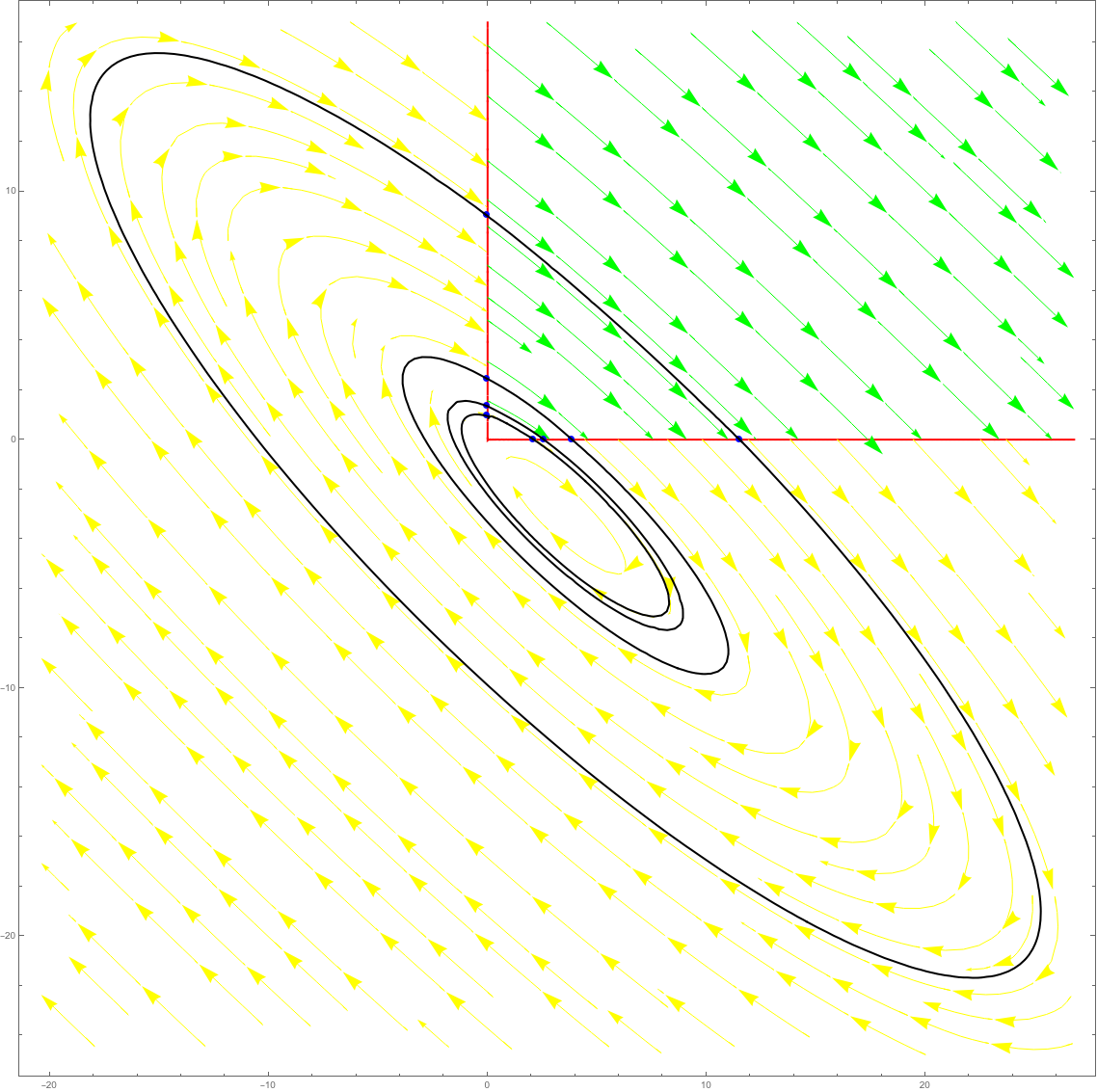}
\caption{The four limit cycles of the discontinuous piecewise differential system \eqref{eq:Xh-61}-\eqref{eq:Yf-61}.}
\label{fig-Hc-H61}
\end{center}
\end{figure}

Now we give an example of a discontinuous piecewise differential system  $(\mathtt{L}_c)$-$(\mathtt{N}_6^{1})$ separated by the nonregular line  $\Sigma$, having four limit cycles. 
In $\Sigma^-$ we  consider the  linear differential center
\begin{equation}\label{eq:Xh-61}
\begin{aligned}
\dot{x} =& \frac{6}{5}+\frac{18}{25}y+2(x+y),\\[4pt]
\dot{y} =& 1-2(x+y),
\end{aligned}
\end{equation}
with the first integral
\begin{equation*}
    H_L(x,y) =2\left(-\frac{x}{2}+\frac{3}{5}y\right)+\frac{9}{25}y^{2}+(x+y)^2.
\end{equation*}
In $\Sigma^{+}$ we consider the Hamiltonian nilpotent saddle of type $(\mathtt{N}_6^{1})$
\begin{equation}\label{eq:Yf-61}
\begin{aligned}
\dot{x} =& \frac{1}{21000}\Big( -218240+36x(9544-5082x+993x^2)+393744y \\
&\quad
+9x(-45808+13989x)y+3(-77840+51789x)y^2+69994y^3\Big), \\[6pt]
\dot{y} =& \frac{1}{7000}\Big( 63232-6x(16832+9x(-1008+193x)) \\
&\quad
-12(9544+3x(-3388+993x))y+3(22904-13989x)y^2-17263y^3\Big),
\end{aligned}
\end{equation}
with the first integral    
\begin{equation*}
\begin{aligned}
 H_{6}^{1}(x,y) =&
\frac{1}{4}\Big(
-\left(-\frac{4}{5}+\frac{3}{10}x+\frac{1}{10}y\right)^4
-\left(-\frac{4}{5}+\frac{3}{5}x+\frac{9}{10}y\right)^4 \Big) \\[4pt]
&-\frac{1}{5}\left(-\frac{4}{5}+\frac{3}{5}x+\frac{9}{10}y\right)^2
-\frac{27}{20}
\left(-\frac{4}{5}+\frac{3}{10}x+\frac{1}{10}y\right)^2 \\
& \left(-\frac{4}{5}+\frac{3}{5}x+\frac{9}{10}y\right)^2 .
\end{aligned}
\end{equation*} 
Solving the system~\eqref{sist_H_Hi}, we obtain four real pairs of solutions
$(p_k,q_k)$, $k=1,\dots,4$, where $p_k=(x_k,0)$ and $q_k=(0,y_k)$. More precisely,
\[
\begin{aligned}
p_1 &\approx (11.5231,0), & q_1 &\approx (0,9.01165),
& p_2 &\approx (3.8652,0), & q_2 &\approx (0,2.44633),\\
p_3 &\approx (2.57797,0), & q_3 &\approx (0,1.3437),
& p_4 &\approx (2.11393,0), & q_4 &\approx (0,0.946661).
\end{aligned}
\]
that provide the crossing limit cycles of discontinuous piecewise differential system \eqref{eq:Xh-61}–\eqref{eq:Yf-61} shown in Figure \ref{fig-Hc-H61}. 

 
 \subsection*{Limit cycles for the discontinuous piecewise differential system $(\mathtt{L}_c)$-$(\mathtt{N}_6^{2})$}
\begin{figure}[t]
\begin{center}
\includegraphics[width=0.45\textwidth, 
height=0.3\textheight]{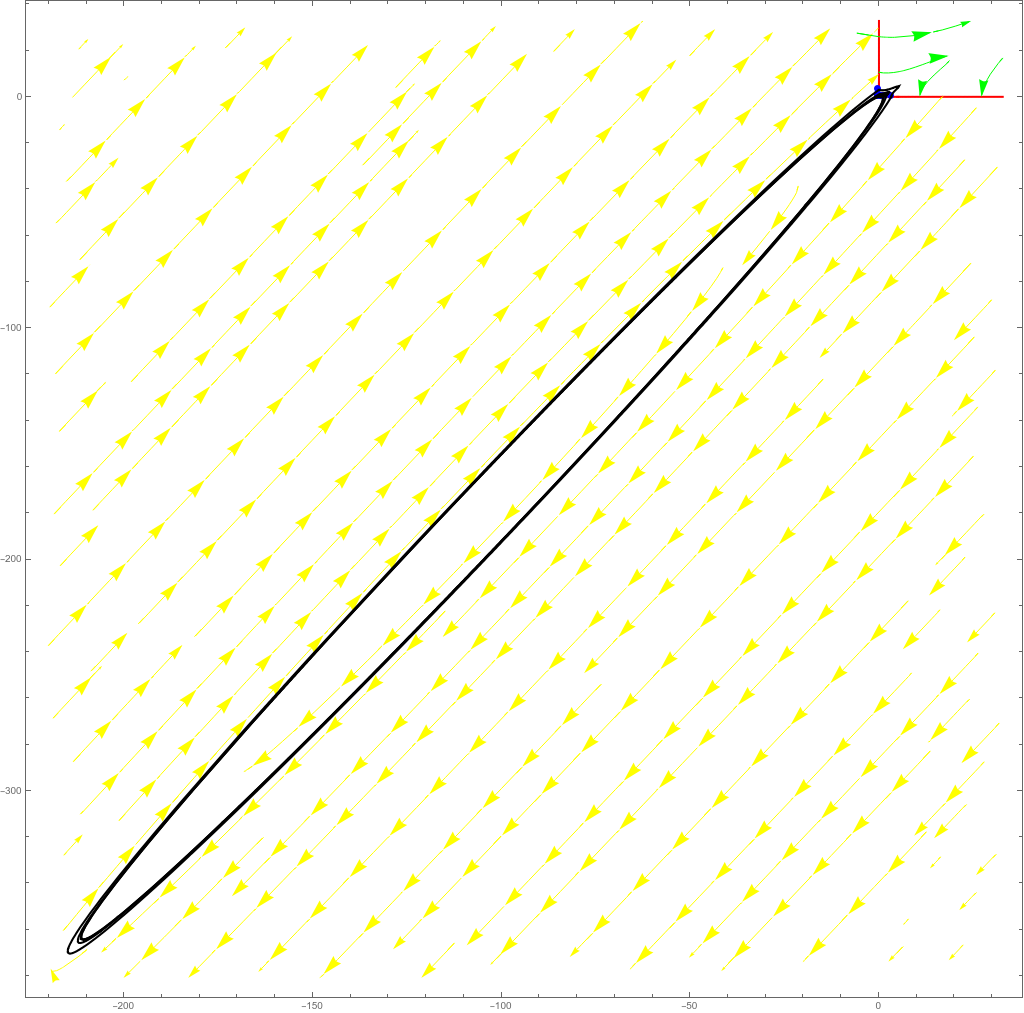}
\includegraphics[width=0.45\textwidth,
height=0.3\textheight]{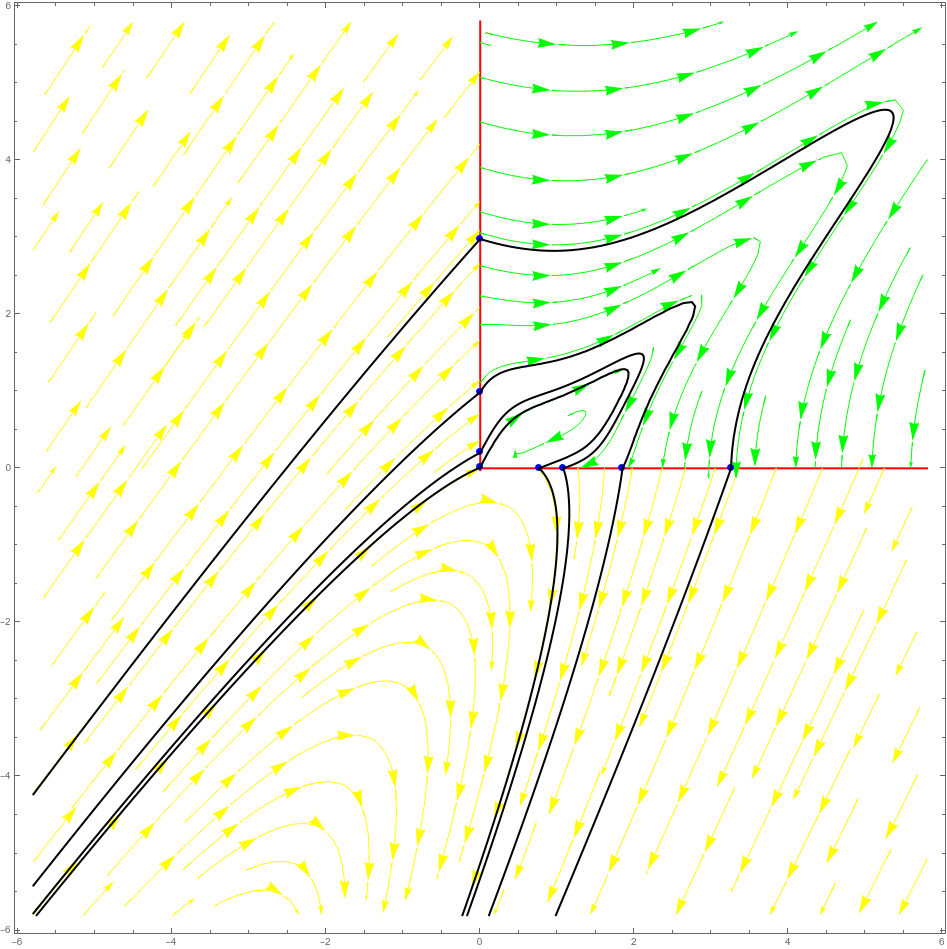}
\caption{The four limit cycles of the discontinuous piecewise differential system \eqref{eq:Xh-62}–\eqref{eq:Yf-62}. Right: zoom around the origin.}
\label{fig-Hc-H62}
\end{center}
\end{figure}

Now we give an example of a discontinuous piecewise differential system  $(\mathtt{L}_c)$-$(\mathtt{N}_6^{2})$ separated by the nonregular line  $\Sigma$, having four limit cycles. 
In $\Sigma^-$ we  consider the  linear differential center
\begin{equation}\label{eq:Xh-62}
\begin{aligned}
\dot{x} =& \frac{87}{50}
-\frac{29}{25}\left(x-\frac{29}{50}y\right)
+\frac{9}{1250}y,\\[4pt]
\dot{y} =& \frac{37}{50}
-2\left(x-\frac{29}{50}y\right),
\end{aligned}
\end{equation}
with the first integral
\begin{equation*}
    H_L(x,y) =\left(x-\frac{29}{50}y\right)^2+2\left(-\frac{37}{100}x+\frac{87}{100}y\right)
+\frac{9}{2500}y^2.
\end{equation*}
In $\Sigma^{+}$ we consider the Hamiltonian nilpotent saddle of type $(\mathtt{N}_6^{2})$

\begin{equation}\label{eq:Yf-62} 
\begin{aligned}
\dot{x} =& \frac{1}{11958250000}\Big( 898774048-5365417169\,x-6417x^{2}\big(-421623+62735x\big) \\
& +3915993356\,y +1302\big(3990158-2625497x\big)xy+1116\big(-5047468 \\
& +987185x\big) y^{2}+3286262384\,y^{3} \Big), \\[6pt]
\dot{y} =& \frac{1}{23916500000}\Big( x\Big(-18163986613+31\big(678017913-261411193x\big)x\Big) \\
& +38502x\big(-281082+62735x\big)y+2604\big(-1995079+2625497x\big)y^{2} \\
& -734465640\,y^{3}+7\big(532650799+1532976334y\big) \Big),
\end{aligned}
\end{equation}
with the first integral

\[
\begin{aligned}
H_{6}^{2}(x,y)=\;&
\frac14\Big(
-\left(\frac{29}{100}-\frac{x}{25}-\frac{27}{50}y\right)^4
-\left(-\frac{51}{100}+\frac{57}{100}x-\frac{1}{50}y\right)^4
\Big) \\[4pt]
& -\frac{31}{200}\left(\frac{29}{100}-\frac{x}{25}-\frac{27}{50}y\right)^2
-\frac{4}{25}
\left(\frac{29}{100}-\frac{x}{25}-\frac{27}{50}y\right) \\[4pt]
& \left(-\frac{51}{100}+\frac{57}{100}x-\frac{1}{50}y\right) -\frac{32}{775}
\left(-\frac{51}{100}+\frac{57}{100}x-\frac{1}{50}y\right)^2 \\[4pt]
& +\frac{12}{25}
\left(\frac{29}{100}-\frac{x}{25}-\frac{27}{50}y\right)^2
\left(-\frac{51}{100}+\frac{57}{100}x-\frac{1}{50}y\right)^2 .
\end{aligned}
\]
Solving the system~\eqref{sist_H_Hi}, we obtain four real pairs of solutions
$(p_k,q_k)$, $k=1,\dots,4$, where $p_k=(x_k,0)$ and $q_k=(0,y_k)$. More precisely,
\[
\begin{aligned}
p_1 &\approx (3.25582,0), & q_1 &\approx (0,2.97642),
& p_2 &\approx (1.84751,0), & q_2 &\approx (0,0.985982),\\
p_3 &\approx (1.07893,0), & q_3 &\approx (0,0.202174),
& p_4 &\approx (0.765476,0), & q_4 &\approx (0,0.0111834).
\end{aligned}
\]
that provide the crossing limit cycles of discontinuous piecewise differential system \eqref{eq:Xh-62}–\eqref{eq:Yf-62} shown in Figure \ref{fig-Hc-H62}. 



\section*{Acknowlegements}

The first author acknowledges partial support from the 2026 Summer Postdoctoral Program at the Instituto de Matemática Pura e Aplicada (IMPA) and from a scholarship granted by the Fundação Arthur Bernardes (FUNARBE).
The second author acknowledges partial support from CNPq under grant No.~169201/2023-6.

\appendix
\section{Crossing polynomials}
\label{app1}
This appendix collects the polynomial expressions arising from the first
integrals in the crossing equations for the six classes of systems studied in
the paper.

\textit{Polynomial $P_1$:}
\begin{align*}
P_1(x,y)=&\frac{1}{4b}\bigl(2a x^{2}-b x^{4}-2b^{2}y^{2}\bigr)
\end{align*}

\textit{Polynomial $P_2$:}
\begin{equation*}
\begin{aligned}
P_2(x,y)=&\frac{1}{2b}\big(
2a^{2}a_{1}c_{1}x
+a^{2}a_{1}^{2}x^{2}
-2a^{2}b_{1}c_{1}y
-a^{2}b_{1}^{2}y^{2} 
+2ab c_{1}x\alpha_{1}
-2b c_{1}^{3}x\alpha_{1} \\
& \quad 
+2aa_{1}bx^{2}\alpha_{1} 
-6a_{1}b c_{1}^{2}x^{2}\alpha_{1}-6a_{1}^{2}b c_{1}x^{3}\alpha_{1}
-2a_{1}^{3}b x^{4}\alpha_{1}
+ b^{2}x^{2}\alpha_{1}^{2}
\\
& \quad 
-2ab c_{1}y\beta_{1}
+2b c_{1}^{3}y\beta_{1}
-2ab b_{1}y^{2}\beta_{1}
+6b b_{1}c_{1}^{2}y^{2}\beta_{1}+6b b_{1}^{2}c_{1}y^{3}\beta_{1} \\
& \quad 
+2b b_{1}^{3}y^{4}\beta_{1}
- b^{2}y^{2}\beta_{1}^{2} +2aa_{1}bx\gamma_{1}
-6a_{1}b c_{1}^{2}x\gamma_{1}  
-6a_{1}^{2}b c_{1}x^{2}\gamma_{1} \\
& \quad 
-2a_{1}^{3}b x^{3}\gamma_{1}-2ab b_{1}y\gamma_{1}
+6b b_{1}c_{1}^{2}y\gamma_{1}
+6b b_{1}^{2}c_{1}y^{2}\gamma_{1}
+2b b_{1}^{3}y^{3}\gamma_{1} \\
& \quad +2b^{2}x\alpha_{1}\gamma_{1}
-2b^{2}y\beta_{1}\gamma_{1}
\big)
\end{aligned}
\end{equation*}

\textit{Polynomial $P_3$:}
\begin{equation*}
\begin{aligned}
P_{3}(x,y)=&\frac{1}{4}\big(
-4a_{1}c c_{1}x
-2a_{1}^{2}c x^{2}
+4b_{1}c c_{1}y
+2b_{1}^{2}c y^{2} -6c_{1}^{2}x^{2}\alpha_{1}^{2}
-12a_{1}c_{1}x^{3}\alpha_{1}^{2} \\
& \quad
-6a_{1}^{2}x^{4}\alpha_{1}^{2}
+ x^{4}\alpha_{1}^{4} +6c_{1}^{2}y^{2}\beta_{1}^{2}
+12b_{1}c_{1}y^{3}\beta_{1}^{2}
+6b_{1}^{2}y^{4}\beta_{1}^{2}
- y^{4}\beta_{1}^{4}
\\
& \quad
-12c_{1}^{2}x\alpha_{1}\gamma_{1}
-24a_{1}c_{1}x^{2}\alpha_{1}\gamma_{1} 
-12a_{1}^{2}x^{3}\alpha_{1}\gamma_{1}
+4x^{3}\alpha_{1}^{3}\gamma_{1}
\\
& \quad
+12c_{1}^{2}y\beta_{1}\gamma_{1}
+24b_{1}c_{1}y^{2}\beta_{1}\gamma_{1}
+12b_{1}^{2}y^{3}\beta_{1}\gamma_{1}
-4y^{3}\beta_{1}^{3}\gamma_{1} 
\\
& \quad
-12a_{1}c_{1}x\gamma_{1}^{2}
-6a_{1}^{2}x^{2}\gamma_{1}^{2}
+12b_{1}c_{1}y\gamma_{1}^{2}
+6b_{1}^{2}y^{2}\gamma_{1}^{2}+6x^{2}\alpha_{1}^{2}\gamma_{1}^{2} \\
& \quad
-6y^{2}\beta_{1}^{2}\gamma_{1}^{2}
+4x\alpha_{1}\gamma_{1}^{3}
-4y\beta_{1}\gamma_{1}^{3}
\big)
\end{aligned}
\end{equation*}

\textit{Polynomial $P_4$:}

\begin{equation*}
\begin{aligned}
P_4(x,y)=&\frac{1}{4b}\big(
4a^{2}a_{1}c_{1}x
+2a^{2}a_{1}^{2}x^{2}
-4a^{2}b_{1}c_{1}y
-2a^{2}b_{1}^{2}y^{2} +4ab c_{1}x\alpha_{1} 
\\
& \quad
+4aa_{1}bx^{2}\alpha_{1}
+2b^{2}x^{2}\alpha_{1}^{2}
-6bc_{1}^{2}x^{2}\alpha_{1}^{2} -12a_{1}bc_{1}x^{3}\alpha_{1}^{2}
-6a_{1}^{2}bx^{4}\alpha_{1}^{2} 
\\
& \quad
+bx^{4}\alpha_{1}^{4} -4ab c_{1}y\beta_{1}
-4abb_{1}y^{2}\beta_{1}
-2b^{2}y^{2}\beta_{1}^{2}
+6bc_{1}^{2}y^{2}\beta_{1}^{2} \\
& \quad
+12bb_{1}c_{1}y^{3}\beta_{1}^{2}
+6bb_{1}^{2}y^{4}\beta_{1}^{2}
-by^{4}\beta_{1}^{4}+4aa_{1}bx\gamma_{1}
-4abb_{1}y\gamma_{1} \\
& \quad
+4b^{2}x\alpha_{1}\gamma_{1} 
-12bc_{1}^{2}x\alpha_{1}\gamma_{1}-24a_{1}bc_{1}x^{2}\alpha_{1}\gamma_{1}
-12a_{1}^{2}bx^{3}\alpha_{1}\gamma_{1}\\
& \quad
+4bx^{3}\alpha_{1}^{3}\gamma_{1} -4b^{2}y\beta_{1}\gamma_{1}
+12bc_{1}^{2}y\beta_{1}\gamma_{1}
+24bb_{1}c_{1}y^{2}\beta_{1}\gamma_{1} \\
& \quad
+12bb_{1}^{2}y^{3}\beta_{1}\gamma_{1}-4by^{3}\beta_{1}^{3}\gamma_{1}
-12a_{1}bc_{1}x\gamma_{1}^{2}
-6a_{1}^{2}bx^{2}\gamma_{1}^{2} \\
& \quad
+12bb_{1}c_{1}y\gamma_{1}^{2} +6bb_{1}^{2}y^{2}\gamma_{1}^{2}
+6bx^{2}\alpha_{1}^{2}\gamma_{1}^{2}
-6by^{2}\beta_{1}^{2}\gamma_{1}^{2}
+4bx\alpha_{1}\gamma_{1}^{3}
\\
& \quad
-4by\beta_{1}\gamma_{1}^{3}
\big)
\end{aligned}
\end{equation*}

\textit{Polynomial $P_5$:}

\begin{equation*}
\begin{aligned}
P_5(x,y)=&\frac{1}{4}\big(
-4 a_1 c c_1 x
-2 a_1^{2} c x^{2}
+4 b_1 c c_1 y
+2 b_1^{2} c y^{2} -6 c_1^{2} x^{2} \alpha_1^{2} \\
& \quad
-12 a_1 c_1 x^{3} \alpha_1^{2}
-6 a_1^{2} x^{4} \alpha_1^{2} 
- x^{4} \alpha_1^{4}+6 c_1^{2} y^{2} \beta_1^{2}
+12 b_1 c_1 y^{3} \beta_1^{2} \\
& \quad
+6 b_1^{2} y^{4} \beta_1^{2}
+ y^{4} \beta_1^{4} -12 c_1^{2} x \alpha_1 \gamma_1
-24 a_1 c_1 x^{2} \alpha_1 \gamma_1
-12 a_1^{2} x^{3} \alpha_1 \gamma_1 \\
& \quad
-4 x^{3} \alpha_1^{3} \gamma_1+12 c_1^{2} y \beta_1 \gamma_1
+24 b_1 c_1 y^{2} \beta_1 \gamma_1
+12 b_1^{2} y^{3} \beta_1 \gamma_1
+4 y^{3} \beta_1^{3} \gamma_1 \\
& \quad -12 a_1 c_1 x \gamma_1^{2}
-6 a_1^{2} x^{2} \gamma_1^{2}
+12 b_1 c_1 y \gamma_1^{2}
+6 b_1^{2} y^{2} \gamma_1^{2} -6 x^{2} \alpha_1^{2} \gamma_1^{2}\\
& \quad
+6 y^{2} \beta_1^{2} \gamma_1^{2}
-4 x \alpha_1 \gamma_1^{3}
+4 y \beta_1 \gamma_1^{3}
\big)
\end{aligned}
\end{equation*}

\textit{Polynomial $P_6$:}
\begin{equation*}
\begin{aligned}
P_{6}(x,y)=&\frac{1}{4b}\big(
4 a^{2} a_1 c_1 x
+2 a^{2} a_1^{2} x^{2}
-4 a^{2} b_1 c_1 y
-2 a^{2} b_1^{2} y^{2}+4 a b c_1 x \alpha_1 \\
& \quad
+4 a a_1 b x^{2} \alpha_1
+2 b^{2} x^{2} \alpha_1^{2}
-6 b c_1^{2} x^{2} \alpha_1^{2} -12 a_1 b c_1 x^{3} \alpha_1^{2}
-6 a_1^{2} b x^{4} \alpha_1^{2} 
\\
& \quad
- b x^{4} \alpha_1^{4}-4 a b c_1 y \beta_1
-4 a b b_1 y^{2} \beta_1
-2 b^{2} y^{2} \beta_1^{2}
+6 b c_1^{2} y^{2} \beta_1^{2} \\
& \quad +12 b b_1 c_1 y^{3} \beta_1^{2}
+6 b b_1^{2} y^{4} \beta_1^{2}
+ b y^{4} \beta_1^{4}+4 a a_1 b x \gamma_1
-4 a b b_1 y \gamma_1 \\
& \quad
+4 b^{2} x \alpha_1 \gamma_1 
-12 b c_1^{2} x \alpha_1 \gamma_1-24 a_1 b c_1 x^{2} \alpha_1 \gamma_1
-12 a_1^{2} b x^{3} \alpha_1 \gamma_1
\\
& \quad
-4 b x^{3} \alpha_1^{3} \gamma_1 -4 b^{2} y \beta_1 \gamma_1
+12 b c_1^{2} y \beta_1 \gamma_1
+24 b b_1 c_1 y^{2} \beta_1 \gamma_1 
\\
& \quad +12 b b_1^{2} y^{3} \beta_1 \gamma_1
+4 b y^{3} \beta_1^{3} \gamma_1 -12 a_1 b c_1 x \gamma_1^{2}
-6 a_1^{2} b x^{2} \gamma_1^{2}
\\
& \quad
+12 b b_1 c_1 y \gamma_1^{2} 
+6 b b_1^{2} y^{2} \gamma_1^{2}-6 b x^{2} \alpha_1^{2} \gamma_1^{2}
+6 b y^{2} \beta_1^{2} \gamma_1^{2}
-4 b x \alpha_1 \gamma_1^{3} \\
& \quad
+4 b y \beta_1 \gamma_1^{3}
\big)
\end{aligned}
\end{equation*}
\textit{Polynomial $P_7$:}
\begin{equation*}
\begin{aligned}
P_7(x,y)=&\frac{1}{4}\big(
-4 a_1 c c_1 x
-4 a_1 c_1^{3} x
-2 a_1^{2} c x^{2}
-6 a_1^{2} c_1^{2} x^{2}
-4 a_1^{3} c_1 x^{3}
- a_1^{4} x^{4} \\
& \quad
+4 b_1 c c_1 y
+4 b_1 c_1^{3} y
+2 b_1^{2} c y^{2}
+6 b_1^{2} c_1^{2} y^{2}
+4 b_1^{3} c_1 y^{3}
+ b_1^{4} y^{4} + x^{4} \alpha_1^{4} \\
& \quad
- y^{4} \beta_1^{4}
+4 x^{3} \alpha_1^{3} \gamma_1
-4 y^{3} \beta_1^{3} \gamma_1 +6 x^{2} \alpha_1^{2} \gamma_1^{2}
-6 y^{2} \beta_1^{2} \gamma_1^{2}
+4 x \alpha_1 \gamma_1^{3} \\
& \quad
-4 y \beta_1 \gamma_1^{3}-6 c_1^{2} x^{2} \alpha_1^{2} \mu
-12 a_1 c_1 x^{3} \alpha_1^{2} \mu
-6 a_1^{2} x^{4} \alpha_1^{2} \mu +6 c_1^{2} y^{2} \beta_1^{2} \mu \\
& \quad
+12 b_1 c_1 y^{3} \beta_1^{2} \mu
+6 b_1^{2} y^{4} \beta_1^{2} \mu -12 c_1^{2} x \alpha_1 \gamma_1 \mu
-24 a_1 c_1 x^{2} \alpha_1 \gamma_1 \mu
\\
& \quad
-12 a_1^{2} x^{3} \alpha_1 \gamma_1 \mu +12 c_1^{2} y \beta_1 \gamma_1 \mu
+24 b_1 c_1 y^{2} \beta_1 \gamma_1 \mu
+12 b_1^{2} y^{3} \beta_1 \gamma_1 \mu \\
& \quad -12 a_1 c_1 x \gamma_1^{2} \mu
-6 a_1^{2} x^{2} \gamma_1^{2} \mu
+12 b_1 c_1 y \gamma_1^{2} \mu
+6 b_1^{2} y^{2} \gamma_1^{2} \mu
\big)
\end{aligned}
\end{equation*}
\textit{Polynomial $P_8$:}
\begin{equation*}
\begin{aligned}
P_8(x,y)=&\frac{1}{4b}\big(
4 a^{2} a_1 c_1 x
-4 a_1 b c_1^{3} x
+2 a^{2} a_1^{2} x^{2}
-6 a_1^{2} b c_1^{2} x^{2}
-4 a_1^{3} b c_1 x^{3}
- a_1^{4} b x^{4} \\
& \quad
-4 a^{2} b_1 c_1 y
+4 b b_1 c_1^{3} y 
-2 a^{2} b_1^{2} y^{2}
+6 b b_1^{2} c_1^{2} y^{2}
+4 b b_1^{3} c_1 y^{3}
+ b b_1^{4} y^{4} \\
& \quad +4 a b c_1 x \alpha_1
+4 a a_1 b x^{2} \alpha_1
+2 b^{2} x^{2} \alpha_1^{2}
+ b x^{4} \alpha_1^{4}-4 a b c_1 y \beta_1
\\
& \quad
-4 a b b_1 y^{2} \beta_1
-2 b^{2} y^{2} \beta_1^{2}
- b y^{4} \beta_1^{4}+4 a a_1 b x \gamma_1 
-4 a b b_1 y \gamma_1
+4 b^{2} x \alpha_1 \gamma_1 \\
& \quad
+4 b x^{3} \alpha_1^{3} \gamma_1 -4 b^{2} y \beta_1 \gamma_1
-4 b y^{3} \beta_1^{3} \gamma_1
+6 b x^{2} \alpha_1^{2} \gamma_1^{2} 
-6 b y^{2} \beta_1^{2} \gamma_1^{2}
\\
& \quad
+4 b x \alpha_1 \gamma_1^{3}
-4 b y \beta_1 \gamma_1^{3}
-6 b c_1^{2} x^{2} \alpha_1^{2} \mu
-12 a_1 b c_1 x^{3} \alpha_1^{2} \mu -6 a_1^{2} b x^{4} \alpha_1^{2} \mu \\
& \quad
+6 b c_1^{2} y^{2} \beta_1^{2} \mu
+12 b b_1 c_1 y^{3} \beta_1^{2} \mu
+6 b b_1^{2} y^{4} \beta_1^{2} \mu -12 b c_1^{2} x \alpha_1 \gamma_1 \mu \\
& \quad
-24 a_1 b c_1 x^{2} \alpha_1 \gamma_1 \mu
-12 a_1^{2} b x^{3} \alpha_1 \gamma_1 \mu +12 b c_1^{2} y \beta_1 \gamma_1 \mu \\
& \quad
+24 b b_1 c_1 y^{2} \beta_1 \gamma_1 \mu
+12 b b_1^{2} y^{3} \beta_1 \gamma_1 \mu -12 a_1 b c_1 x \gamma_1^{2} \mu
-6 a_1^{2} b x^{2} \gamma_1^{2} \mu \\
& \quad
+12 b b_1 c_1 y \gamma_1^{2} \mu
+6 b b_1^{2} y^{2} \gamma_1^{2} \mu
\big)
\end{aligned}
\end{equation*}
\textit{Polynomial $P_9$:}
\begin{equation*}
\begin{aligned}
P_9(x,y)=&\frac{1}{4}\big(
-4 a_1 c c_1 x
-4 a_1 c_1^{3} x
-2 a_1^{2} c x^{2}
-6 a_1^{2} c_1^{2} x^{2}
-4 a_1^{3} c_1 x^{3}
- a_1^{4} x^{4} \\
& \quad +4 b_1 c c_1 y
+4 b_1 c_1^{3} y
+2 b_1^{2} c y^{2}
+6 b_1^{2} c_1^{2} y^{2}
+4 b_1^{3} c_1 y^{3}
+ b_1^{4} y^{4}
- x^{4} \alpha_1^{4} \\
& \quad
+ y^{4} \beta_1^{4}
-4 x^{3} \alpha_1^{3} \gamma_1
+4 y^{3} \beta_1^{3} \gamma_1 -6 x^{2} \alpha_1^{2} \gamma_1^{2}
+6 y^{2} \beta_1^{2} \gamma_1^{2}
-4 x \alpha_1 \gamma_1^{3} 
\\
& \quad
+4 y \beta_1 \gamma_1^{3}-6 c_1^{2} x^{2} \alpha_1^{2} \mu
-12 a_1 c_1 x^{3} \alpha_1^{2} \mu
-6 a_1^{2} x^{4} \alpha_1^{2} \mu+6 c_1^{2} y^{2} \beta_1^{2} \mu \\
& \quad
+12 b_1 c_1 y^{3} \beta_1^{2} \mu
+6 b_1^{2} y^{4} \beta_1^{2} \mu -12 c_1^{2} x \alpha_1 \gamma_1 \mu
-24 a_1 c_1 x^{2} \alpha_1 \gamma_1 \mu\\
& \quad
-12 a_1^{2} x^{3} \alpha_1 \gamma_1 \mu +12 c_1^{2} y \beta_1 \gamma_1 \mu
+24 b_1 c_1 y^{2} \beta_1 \gamma_1 \mu
+12 b_1^{2} y^{3} \beta_1 \gamma_1 \mu \\
& \quad -12 a_1 c_1 x \gamma_1^{2} \mu
-6 a_1^{2} x^{2} \gamma_1^{2} \mu
+12 b_1 c_1 y \gamma_1^{2} \mu
+6 b_1^{2} y^{2} \gamma_1^{2} \mu
\big)
\end{aligned}
\end{equation*}
\textit{Polynomial $P_{10}$:}
\begin{equation*}
\begin{aligned}
P_{10}(x,y)=&\frac{1}{4b}\big(
4 a^{2} a_1 c_1 x
-4 a_1 b c_1^{3} x
+2 a^{2} a_1^{2} x^{2}
-6 a_1^{2} b c_1^{2} x^{2}
-4 a_1^{3} b c_1 x^{3}
- a_1^{4} b x^{4} \\
& \quad-4 a^{2} b_1 c_1 y
+4 b b_1 c_1^{3} y
-2 a^{2} b_1^{2} y^{2}
+6 b b_1^{2} c_1^{2} y^{2}
+4 b b_1^{3} c_1 y^{3}
+ b b_1^{4} y^{4} \\
& \quad +4 a b c_1 x \alpha_1
+4 a a_1 b x^{2} \alpha_1
+2 b^{2} x^{2} \alpha_1^{2}
- b x^{4} \alpha_1^{4} -4 a b c_1 y \beta_1 \\
& \quad
-4 a b b_1 y^{2} \beta_1 
-2 b^{2} y^{2} \beta_1^{2}
+ b y^{4} \beta_1^{4}+4 a a_1 b x \gamma_1
-4 a b b_1 y \gamma_1 \\
& \quad
+4 b^{2} x \alpha_1 \gamma_1
-4 b x^{3} \alpha_1^{3} \gamma_1 \-4 b^{2} y \beta_1 \gamma_1
+4 b y^{3} \beta_1^{3} \gamma_1
-6 b x^{2} \alpha_1^{2} \gamma_1^{2}
\\
& \quad
+6 b y^{2} \beta_1^{2} \gamma_1^{2}-4 b x \alpha_1 \gamma_1^{3} +4 b y \beta_1 \gamma_1^{3} 
-6 b c_1^{2} x^{2} \alpha_1^{2} \mu
-12 a_1 b c_1 x^{3} \alpha_1^{2} \mu \\
& \quad
-6 a_1^{2} b x^{4} \alpha_1^{2} \mu
+6 b c_1^{2} y^{2} \beta_1^{2} \mu
+12 b b_1 c_1 y^{3} \beta_1^{2} \mu  
+6 b b_1^{2} y^{4} \beta_1^{2} \mu
\\
& \quad
-12 b c_1^{2} x \alpha_1 \gamma_1 \mu 
-24 a_1 b c_1 x^{2} \alpha_1 \gamma_1 \mu
-12 a_1^{2} b x^{3} \alpha_1 \gamma_1 \mu\\
& \quad +12 b c_1^{2} y \beta_1 \gamma_1 \mu
+24 b b_1 c_1 y^{2} \beta_1 \gamma_1 \mu
+12 b b_1^{2} y^{3} \beta_1 \gamma_1 \mu\\
& \quad -12 a_1 b c_1 x \gamma_1^{2} \mu
-6 a_1^{2} b x^{2} \gamma_1^{2} \mu
+12 b b_1 c_1 y \gamma_1^{2} \mu
+6 b b_1^{2} y^{2} \gamma_1^{2} \mu
\big)
\end{aligned}
\end{equation*}


%
\bibliographystyle{acm}
\bibliography{sample}

\end{document}